\documentclass[pdftex,12pt,a4paper,reqno]{amsart}
\usepackage[T1]{fontenc}
\usepackage[utf8]{inputenc}
\usepackage[english]{babel}
\usepackage{lmodern}
\usepackage{amsthm}
\usepackage{amsmath}
\usepackage{amsfonts}
\usepackage{amssymb}
\usepackage[all]{xy}
\usepackage{calrsfs}
\usepackage{slashed}
\usepackage{babel}
\usepackage{array}
\usepackage{tensor}
\theoremstyle{plain}
\newtheorem{theo}{Theorem}[section]
\newtheorem{prop}[theo]{Proposition}
\newtheorem{lemm}[theo]{Lemma}

\newtheorem{cor}[theo]{Corollary}

\theoremstyle{definition}
\newtheorem{defi}[theo]{Definition}

\theoremstyle{remark}

\newtheorem{ex}[theo]{Example}

\newtheorem{rk}[theo]{Remark}

\newcommand{\C}{\mathbb{C}}

\newcommand{\N}{\mathbb{N}}

\newcommand{\R}{\mathbb{R}}

\newcommand{\Sph}{\mathbb{S}}

\usepackage[bookmarks,bookmarksopen,bookmarksdepth=2]{hyperref}
\usepackage[usestackEOL]{stackengine}
\usepackage{etoolbox}
\patchcmd{\thmhead}{(#3)}{#3}{}{}
\usepackage[left=3cm,right=3cm,top=3cm,bottom=3cm]{geometry}
\usepackage{mdframed}

\title{}
\author{}

\everymath{\displaystyle\everymath{}}

\begin{document}

%\maketitle

\setcounter{tocdepth}{2}

%\tableofcontents

%\newpage

\title{CR-harmonic maps}
\author{Gautier Dietrich}
\thanks{The author was supported in part by the grant ANR-17-CE40-0034 of the French National Research Agency ANR (project CCEM)}
\address{Institut Montpelliérain Alexander Grothendieck\\ Université de Montpellier\\ CNRS\\ Case courrier 051\\ Place Eugène Bataillon\\ 34090 Montpellier\\ France} 
\address{Université Paul-Valéry Montpellier 3%\\ UFR 6\\ Route de Mende \\34199 Montpellier cedex 5\\ France
} 
\email{gautier.dietrich@ac-montpellier.fr}
\date{}

\begin{abstract}
%We develop the notion of renormalized energy in CR geometry, as the complex counterpart of the work of V. Bérard in the real case. This energy is a CR invariant functional on maps from a strictly pseudoconvex pseudohermitian compact manifold to a Riemannian manifold. Its critical points, which we call CR-harmonic maps, satisfy a CR covariant PDE, and generalize subharmonic maps. The induced operator coincides on functions with the CR Paneitz operator.
We develop the notion of renormalized energy in CR geometry, for maps from a strictly pseudoconvex pseudohermitian manifold to a Riemannian manifold. %, as the complex counterpart of the work of V. Bérard in the real case
This energy is a CR invariant functional, whose critical points, which we call CR-harmonic maps, satisfy a CR covariant partial differential equation. % and generalize CR-holomorphic maps
The corresponding operator coincides on functions with the CR Paneitz operator.
\end{abstract}

\maketitle

\setcounter{tocdepth}{1}
\tableofcontents

\section{Introduction}

Let $(M,g)$ and $(N,h)$ be two Riemannian manifolds. The \emph{Dirichlet energy} of a map $\varphi:(M,g)\rightarrow(N,h)$ is defined as $$E(\varphi)=\frac{1}{2}\int_M\|T\varphi\|_{g,h}^2 d\mathrm{vol}_g.$$ When $\dim M=2$, the energy is conformally invariant with respect to $g$. This is of considerable usefulness, e.g. to construct conformal minimal immersions of Riemann surfaces \cite{Mil79}. However, in higher dimension, the energy is no longer conformally invariant.

Critical points of a functional are solutions to a partial differential equation called the \emph{Euler-Lagrange equation} of the functional; in other words, they form the kernel of a certain differential operator. In our case, the critical points of the Dirichlet energy are called \emph{harmonic maps}, and harmonic \emph{functions} $\varphi:(M,g)\rightarrow(\R,\mathrm{eucl})$ coincide with the kernel of the Laplacian.

%This notion is defined for a non-compact, \emph{conformally compact} Riemannian manifold $(X,g)$. That means that $X$ is the interior of a compact manifold with boundary $\overline{X}=X\cup\partial\overline{X}$ and that, for all choice of $r\in C^\infty(\overline{X},\R_+)$ defining $M:=\partial\overline{X}$, meaning that $r>0$ on $X$, $r=0$ and $dr\neq 0$ on $\lbrace 0\rbrace\times M$, the metric $r^2 g$ can be regularly extended to a metric on $\overline{X}$. Conformal compactness being independent of the choice of $r$, the boundary $M$ is naturally equipped with the conformal structure $[r^2g|_{TM}]$, called \emph{conformal infinity} of $(X,g)$. If moreover $|dr|_{r^2g}=1$ on $M$, then the sectional curvatures of $g$ go to $-1$ at infinity. The manifold $(X,g)$ is then called \emph{asymptotically hyperbolic}.

%Vice versa, given a conformal manifold $(M,[g_0])$ of dimension $n$, there exists an asymptotically hyperbolic manifold $(X,g)$ with conformal infinity $(M,[g_0])$ which is Einstein up to order $N=\infty$ if $n$ is odd, $N=n-2$ if $n$ is even, that is $$\mathrm{Ric}(g)=-ng+O(r^n),$$ and $g$ is unique modulo $O(r^N)$. $g$ is then called a \emph{Poincaré metric}.

%A \emph{renormalized energy} has been introduced by V. Bérard, using asymptotically hyperbolic geometry \cite{Ber13}. Namely, 
In a recent work, V. Bérard has shown the existence, given two Riemannian manifolds $(M,g)$ and $(N,h)$, with $M$ of even dimension $n$, of a functional $\mathcal{E}_g^n$ on $C^\infty(M,N)$, conformally invariant with respect to $g$, and equal to the usual energy when $n=2$ \cite{Ber13}. This functional is called \emph{renormalized energy}, and its critical points are called \emph{conformal-harmonic maps}. %When $(M,g)$ is an Einstein manifold, c
Conformal-harmonic maps generalize harmonic maps; moreover, when $n=4$ and $N=\R$, the induced operator coincides with the Paneitz operator.

We develop here the notions of CR-harmonicity and renormalized energy in CR geometry. %When the source pseudohermitian manifold is normal and pseudo-Einstein, 
CR-harmonic maps also generalize CR-holomorphic maps, which are notoriously hard to come by. 
When $\dim M=3$ and $N=\R$, the induced operator coincides with the CR Paneitz operator. This generalizes the recent work of T. Marugame \cite{Mar18}. Another extension of the CR Paneitz operator to maps has been proposed by T. Chong, Y. Dong, Y. Ren, and G. Yang \cite{CDRY17}. %On functions, we retrieve the CR GJMS operators of A. R. Gover and C. R. Graham \cite{GG05}. 
The main result is the following, which summarizes Proposition \ref{ene} and Theorem \ref{obs}:
\begin{theo}\label{corpsA}
Let $(M^{2n+1},H,J,\theta)$ be a compact strictly pseudoconvex pseudohermitian manifold and $(N,h)$ be a Riemannian manifold. There exists a functional $F_n$ on $C^\infty(M,N)$ which is a CR invariant, \emph{i.e.} conformally invariant with respect to $\theta$. For $\varphi\in C^\infty(M,N)$, it reads 
\begin{displaymath}
\begin{array}{rcl}
F_n(\varphi) &=& \frac{(-1)^{n+1}}{2n!^2}\int_M \left<(\delta^{\theta,h}_b\nabla^{\varphi^*h})^{n-1}\delta^{\theta,h}_bT\varphi,\delta^{\theta,h}_bT\varphi\right>_h \theta\wedge d\theta^n\\
&&+\text{lower order terms (in derivatives of $\varphi$)},
\end{array}
\end{displaymath}
where $\delta^{\theta,h}_b$ is the Webster divergence on $\Omega^1(M)\otimes\varphi^* TN$.
%\end{theo}
%
%\begin{theo}

The Euler-Lagrange equation of $F_n$ is a partial differential equation of order $2n+2$, itself CR covariant. For $\varphi\in C^\infty(M,N)$, it reads $$0=\frac{(-1)^n}{n!}(\delta^{\theta,h}_b\nabla^{\varphi^*h})^n \delta^{\theta,h}_bT\varphi+\text{lower order terms (in derivatives of $\varphi$)}.$$% that is, the leading term of the corresponding operator is $\frac{(-1)^n}{n!}(\delta^{\theta,h}_b\nabla^{\varphi^*h})^n \delta^{\theta,h}_bT.$
\end{theo}
\noindent Moreover, we provide explicit computations of $P_1$ and $F_1$ in Theorems \ref{corpsB1} and \ref{corpsB2} respectively.

The paper is organized as follows: in Section \ref{ache}, we recall notions of asymptotically complex hyperbolic geometry. In Section \ref{sechar}, we adapt the classical construction by C. R. Graham, R. Jenne, L. J. Mason, and G. A. J. Sparling to obtain a CR Paneitz operator acting on maps, and we define CR-harmonicity \cite{GJMS92}. We also provide an explicit computation of the operator in dimension $3$. In Section \ref{secren}, we develop the corresponding notion of renormalized energy. Section \ref{secfur} presents computations in higher dimension, which do not allow for an explicit expression of the operator. Finally, Section \ref{secfef} gives a correspondence between CR-harmonic maps on a pseudohermitian manifold and conformal-harmonic maps on its Fefferman bundle.

We adopt the following convention: small Greek letters will denote indices in $\lbrace 1,\ldots,n\rbrace$; capital Greek letters, in $\lbrace 1,\ldots,n,\overline{1},\ldots,\overline{n}\rbrace$; small Latin letters, in $\lbrace 0,1,\ldots,n\rbrace$; capital Latin letters, in $\lbrace 0,1,\ldots,n,\overline{0},\overline{1},\ldots,\overline{n}\rbrace$. Moreover, we use the Einstein summation convention everywhere.

\textbf{Acknowledgements.} I am deeply grateful to my supervisor, Marc Herzlich, for introducing me to these questions, for his numerous advices and his precious help. I also profusely thank the anonymous referee of this paper and the referees of my PhD thesis, Jih-Hsin Cheng and Colin Guillarmou, for their careful reading.

\section{ACHE manifolds}\label{ache}

\emph{Asymptotically hyperbolic} manifolds (AH for short) are manifolds which admit a \emph{conformal infinity}, that is to say a boundary equipped with a conformal structure which is, roughly speaking, a generalization of the standard conformal sphere seen as the boundary of the Poincaré disk. Reciprocally, every compact conformal manifold can be filled with an AH manifold $X^{n+1}$ whose metric is Einstein, thus called \emph{AH-Einstein} or \emph{AHE}, when $n$ is odd. When $n$ is even, a conformally invariant obstruction to the existence of a smooth up to the boundary AHE metric appears \cite{FG85,GH05}. Recently, M. J. Gursky and G. Székelyhidi have announced that an AHE metric exists locally for all $n\geq 3$ \cite{GS17}. % The link between conformal and asymptotically hyperbolic geometry is fruitful, particularly when looking for conformal invariants. It is also central in particle physics, where it is known as \emph{AdS/CFT correspondence}.
This approach provides a correspondence between a Riemannian structure on a manifold and a conformal structure on its boundary. Information on the conformal infinity can thus be read on the AHE metric.

The complex counterparts of AH manifolds, asymptotically \emph{complex} hyperbolic manifolds  (ACH for short), have been introduced by C. Epstein, R. Melrose, and G. Mendoza \cite{EMM91}. They generalize the construction by C. Fefferman, S.-Y. Cheng, and S.-T. Yau, of \emph{asymptotically Bergman metrics}, which are Kähler-Einstein metrics on bounded strictly pseudoconvex domains of $\C^{n+1}$, which are asymptotic to the CR structure of the boundary %,  
\cite{Fef76,CY80}. The regularity of these metrics near the boundary has been studied by J. Lee and R. Melrose \cite{LM82}. To an ACH manifold thus corresponds a \emph{CR infinity}. For example, the CR infinity of the complex hyperbolic space $\C\bold{H}^{n+1}$ is $\Sph^{2n+1}$ endowed with its standard CR structure.

Because of the anisotropy of their structure, pseudohermitian manifolds of odd dimension $N$ often behave, \emph{mutatis mutandis}, like Riemannian manifolds of dimension $N+1$. They are sometimes said to have \emph{homogeneous dimension} $N+1$ \cite{JL89}. In particular, ACH manifolds have been known to share similarities with the "$n$ even" real case. The asymptotic development of ACH-Einstein and -Kähler-Einstein metrics has been extensively studied by O. Biquard, M. Herzlich, and Y. Matsumoto, and obstructions to smoothness have been identified \cite{Biq00, BH05, Mat14}.\\

Let us consider the sphere $\Sph^{2n+1}\subset\C^{n+1}$ endowed with its standard contact form $$\theta_0=\frac{i}{4}\left(z_jd\overline{z}^j-\overline{z}_jdz^j\right)|_{\Sph^{2n+1}}.$$ Let $\gamma_0=d\theta_0(\cdot,i\cdot)$ be the induced metric on the contact distribution $\ker\theta_0$. The \emph{Bergman metric} on the ball $\mathbb{B}^{2n+2}$ is given in polar coordinates by $$g_0=dt^2+4\sinh^2(t)\theta_0^2+4\sinh^2\left(\frac{t}{2}\right)\gamma_0.$$
This metric is Kähler and has constant holomorphic sectional curvature $-1$. The space $(\mathbb{B}^{2n+2}, g_0)$ is known as the \emph{complex hyperbolic space} and is denoted by $\C\textbf{H}^{n+1}$.

More generally, let $(M,H,J)$ be a $(2n+1)$-dimensional orientable compact strictly pseudoconvex CR manifold. Namely, $H$ is an orientable hyperplane distribution in $TM$ and $J$ is a complex structure on $H$. Let $\theta$ be a compatible positive contact form and $\gamma = d\theta(\cdot,J\cdot)$ be the induced metric. Let $R$ be the Reeb field. Let $\nabla^\theta$ be the Tanaka-Webster connection of $(M,H,J,\theta)$ and $\tau$ be the pseudohermitian torsion.
% $\tau:=T^W(R,\cdot)$.

%$\gamma^\tau:=i\gamma_{\alpha\overline{\beta}}\left(\tau^{\overline{\beta}}\wedge\theta^\alpha-\tau^\alpha\wedge\theta^{\overline{\beta}}\right)=i\left(\tau_{\alpha\beta}\theta^\alpha\wedge\theta^\beta-\tau_{\overline{\alpha}\overline{\beta}}\theta^{\overline{\alpha}}\wedge\theta^{\overline{\beta}}\right)$ and 

Let %$X$ be a $(2n+2)$-dimensional manifold such that the complement of some compact set is diffeomorphic to 
$\overline{X}=[0,\varepsilon)\times M$, let $\pi:\overline{X}\rightarrow M$ be the natural projection, and let $r$ be the coordinate on $[0,\varepsilon)$. Let $X$ be the interior of $\overline{X}$. Let $g_0$ be the metric on $X$ $$g_0 = \frac{dr^2}{r^2}+\frac{\theta^2}{r^2}+\frac{\gamma}{r}.$$

A function $s\in C^\infty(\overline{X},\R_+)$ is called \emph{boundary defining} if $s>0$ on $X$, $s=0$ and $ds\neq 0$ on $\lbrace 0\rbrace\times M$. Equivalently, $s=e^{f} r$ for some $f$ in $C^\infty(\overline{X},\R)$. A conformal change of the boundary defining function corresponds to a conformal change of the contact form. Indeed, let us consider $g_0$ as $g_0(r,\theta)$, then, for $f$ in $C^\infty(\overline{X},\R)$, $$g_0(e^{f} r,\theta)=g_0(r,e^{-f\vert_M}\theta).$$

We define an order $O_e$ adapted to $g_0$. A normal basis with respect to $g_0$ is $e=\left(r\partial_r,rR,r^\frac{1}{2}T_A\right)$, where $\left(T_A\right)$ is an orthonormal basis for $\gamma$, considered as a Hermitian metric. Its dual basis is $e^*=\left(r^{-1}dr,r^{-1}\theta,r^{-\frac{1}{2}}\theta^\alpha,r^{-\frac{1}{2}}\theta^{\overline{\alpha}}\right)$. The order $O_e$ takes $e$ and $e^*$ for reference. Thus, we have for example $$\gamma=\theta^\alpha\circ\theta^{\overline{\alpha}}=r \left(r^{-\frac{1}{2}}\theta^\alpha\right)\circ\left(r^{-\frac{1}{2}}\theta^{\overline{\alpha}}\right)=O_e(r),$$ where $\lambda\circ\mu:=\lambda\otimes\mu+\mu\otimes\lambda$.

%Let $C_\delta^\infty$ be the space of smooth functions $u$ on $X$ such that $r^{-\delta}\nabla^k u$ is bounded for any $k$.

\begin{defi}[\cite{Biq00}]
A metric $g$ on $X$ is called \emph{asymptotically complex hyperbolic}, or \emph{ACH}, if $g-g_0=o_e(1)$. %belongs to $C_\delta^\infty$ for some $\delta>0$. 
The CR manifold $(M,H,J)$ is then called the \emph{CR infinity} of $(X,g)$.
\end{defi}

\begin{ex}\label{exlam}
For $\lambda>0$, $$g = \frac{dr^2}{r^2}+\frac{(1-\lambda^2r^2)^2}{r^2}\theta^2+\frac{(1-\lambda r)^2}{r}\gamma$$ is an ACH metric on $X$. Moreover, if $(M,H,J,\theta)$ is \emph{Einstein}, \emph{i.e.} pseudo-Einstein with vanishing pseudohermitian torsion, with $\mathrm{Ric}_W(J,\theta)=2(n+1)\lambda\gamma$, then $g$ is an Einstein metric, satisfying $$\mathrm{Ric}(g)=-\frac{n+2}{2}g.$$ Indeed, a complex structure $\tilde{J}$ compatible with $g$ on $X$ is given by $\tilde{J}|_{H\times\lbrace r\rbrace}=J$ and $\tilde{J}\partial_r=-\frac{R}{1-\lambda^2r^2}$, \emph{i.e.} $dr\circ\tilde{J}=(1-\lambda^2r^2)\theta$. %a direct computation (performed with SageManifolds) shows that $$\mathrm{Ric}(g)=-\frac{n+2}{2}g-\frac{n}{2}\frac{(1+\lambda r)^4}{r^2}\theta^2+\frac{(1+\lambda r)^2}{2r}\gamma-2\lambda(n+1)\gamma.$$
Let $\theta^0:=\frac{1}{\sqrt{2}}\left(\frac{1}{1-\lambda^2r^2}dr-i\theta\right)$ and let $\sigma:=\theta^0\wedge\theta^1\wedge\ldots\wedge\theta^n$ be a section of the canonical bundle.
Then $$d\sigma=\frac{i}{\sqrt{2}} d\theta\wedge\theta^1\wedge\ldots\wedge\theta^n-\theta^0\wedge d\theta^1\wedge\ldots\wedge\theta^n+\ldots+(-1)^n\theta^0\wedge \theta^1\wedge\ldots\wedge d\theta^n,$$
where the first term vanishes and, since $\tau=0$, $d\theta^\alpha=\theta^\beta\wedge\omega_\beta^\alpha$, hence
$$d\sigma=-\omega_\alpha^\alpha\wedge\sigma.$$ The curvature form of $\sigma$, in the sense of \cite{BH05}, is hence given by $-d\omega_\alpha^\alpha=-\mathcal{R}\indices{^\theta_\alpha^\rho_{\rho\overline{\alpha}}}\theta^\alpha\wedge\theta^{\overline{\alpha}}=2i(n+1)\lambda d\theta$. Moreover, $$\sigma\wedge\overline{\sigma}=\frac{(-1)^{n+1}ir^{n+2}}{(1-\lambda^2r^2)^2(1-\lambda r)^{2n}}\left(r^{-1}dr\right)\wedge\left((1-\lambda^2r^2)r^{-1}\theta\right)\wedge\left((1-\lambda r)r^{-\frac{1}{2}}\theta^1\right)\wedge\ldots\wedge\left((1-\lambda r)r^{-\frac{1}{2}}\theta^{\overline{n}}\right).$$ Consequently, $$|\sigma|_g^2=\frac{r^{n+2}}{(1-\lambda^2r^2)^2(1-\lambda r)^{2n}},$$
hence $\ln|\sigma|_g^2=(n+2)\ln r-2\ln(1+\lambda r)-(2n+2)\ln(1-\lambda r)$.
We have $$\partial r=\frac{1}{2}\left(dr-i(1-\lambda^2r^2)\theta\right),$$ hence $$i\overline{\partial}\partial r=-\lambda^2r dr\wedge\theta+\frac{1-\lambda^2r^2}{2}d\theta\quad\text{and}\quad i\overline{\partial} r\wedge\partial r=\frac{1-\lambda^2r^2}{2}dr\wedge\theta.$$
The Ricci form of $g$ is then given by 
\begin{displaymath}
\begin{array}{rcl}
\rho_g&=&-i\partial\overline{\partial}\ln|\sigma|_g^2+id\omega_\alpha^\alpha\\
&=&-i\partial\overline{\partial}\ln|\sigma|_g^2+2(n+1)\lambda d\theta\\
&=&\frac{n+2}{2}\left(\frac{1-\lambda^2r^2}{r^2}dr\wedge\theta-\frac{(1-\lambda r)^2}{r}d\theta\right).
\end{array}
\end{displaymath}
\end{ex}

With this example in mind, one may ask if there is in general an ACH \emph{Einstein} (ACHE for short) metric on $X$. Contrarily to the theorem of Cheng-Yau for domains of $\C^{n+1}$, such a metric may not exist in general \cite{CY80}. Nevertheless, there are \emph{formally determined} almost ACHE metrics, in the following sense:

\begin{defi}
In any asymptotic development $\sum_k a_k(p) r^k$, the term $a_k$, seen as a function on $M$, is called \emph{formally determined} if it is a universal polynomial on a finite jet of the CR structure at $p\in M$ only.
\end{defi}

\begin{theo}[\cite{Mat14}]\label{poinc}
There is an ACH metric $g_E$ on $X$, which is Einstein up to order $n+1$, i.e. $$\mathrm{Ric}(g_E)=-\frac{n+2}{2}g_E+O_e(r^{n+1}),$$ where $O_e$ denotes the order with respect to any basis $e$ orthonormal for $g_0$. The metric $g_E$ is formally determined modulo $O_e(r^{n+1})$. % and of the form $$g=\frac{dr^2+g_r}{r^2}.$$
Moreover, we have the asymptotic development
$$g_E = g_0 + \Phi + O_e(r^\frac{3}{2}),$$ where $$\Phi = -2 \mathrm{Sch}_W(J,\theta) + 2 \gamma(J\tau\cdot,\cdot),$$ where $$\mathrm{Sch}_W(J,\theta)=\frac{1}{n+2}\left(\mathrm{Ric}_W(J,\theta)-\frac{\mathrm{Scal}_W(J,\theta)}{2(n+1)}\gamma\right)$$ is the \emph{CR Schouten tensor}.
\end{theo}
\begin{rk} Note that $\Phi=O_e(r)$. \end{rk}

We thus have a formally determined \emph{almost} ACHE metric on $X$. A more convenient metric for our study would be an almost ACH-\emph{Kähler}-Einstein metric on $X$. We have at hand the following results:

\begin{prop}[\cite{BH05}]\label{comp2}
One can construct on $X$ a formal complex structure $J_X$, entirely formally determined by the CR infinity, starting from the almost complex structure $\tilde{J}$, which is the extension of $J$ to $X$ with $\tilde{J}\partial_r=R$. Moreover, an extension $\tilde{\nabla}^\theta$ of $\nabla^\theta$ to $X$ is given by $$\tilde{\nabla}^\theta r\partial_r=\tilde{\nabla}^\theta rR=\tilde{\nabla}^\theta_{r\partial_r}r^\frac{1}{2}T_A=0.$$ Let $\tilde{T}^\theta$ be the torsion of $\tilde{\nabla}^\theta$ and $\tilde{\tau}:=\iota_R \tilde{T}^\theta$. An asymptotic development of $J_X$ is then given by $$J_X=\tilde{J}-2r\tilde{\tau}+O_e(r^\frac{5}{2}).$$%$$J_X=\tilde{J}-2r\tau-r^2 J\tilde{\nabla}^\theta_R\tau +2r^2|\tau|^2+O_e(r^\frac{7}{2}).$$
\end{prop}

\begin{theo}[\cite{Fef76,BH05, Her07}]\label{ke}
There is a formally determined ACH Kähler metric $g_{KE}$ on $(X,J_X)$, which is Einstein up to order $n+\frac{3}{2}$, i.e. $$\mathrm{Ric}(g_{KE})=-\frac{n+2}{2}g_{KE}+O_e(r^{n+\frac{3}{2}}).$$ Moreover, $g_E$ and $g_{KE}$ coincide up to order $n+\frac{1}{2}$.
\end{theo}

In dimension $2n+1=3$, the asymptotic development of $g_{KE}$, and therefore of $g_E$, is known at order $\frac{3}{2}$, which will be essential in Sections \ref{exob} and \ref{exen}: %for the explicit characterisation of CR harmonic maps:

\begin{theo}[\cite{BH05,Her07}]
\label{expl}
When $n=1$, we have the asymptotic development $$g_{KE}=g_0+\Phi_{AB}\theta^A\circ\theta^B+\Psi_{0\overline{1}}\theta^0\circ\theta^{\overline{1}}+\Psi_{\overline{0}1}\theta^{\overline{0}}\circ\theta^1+O_e(r^2),$$ where $$\Psi_{0\overline{1}}=-\sqrt{2}\left(\frac{1}{6}{\mathrm{Scal}_W}_{,\overline{1}}-\frac{2i}{3}\tau^1_{\overline{1},1}\right),$$ and $\Phi$ is given by Theorem \ref{poinc}: $$\Phi_{1\overline{1}}=-\frac{\mathrm{Scal}_W}{4} \quad \text{and}\quad \Phi_{11}=-i \tau_1^{\overline{1}}.$$
\end{theo}

\section{CR-harmonic maps}\label{sechar}

\subsection{Definitions}

Let $(M,H,J)$ be a $(2n+1)$-dimensional orientable, compact, strictly pseudoconvex CR manifold and $(X,g)$ be an ACH manifold with CR infinity $(M,H,J)$, where $g$ is the approximately ACH-Kähler-Einstein metric given by Theorem \ref{ke}. Let $\pi:X\rightarrow M$ be the standard projection. Let $(N,h)$ be a Riemannian manifold. Let $\varphi \in C^\infty(M,N)$, and let $\tilde{\varphi} \in C^\infty(\overline{X},N)$ be any \emph{extension} of $\varphi$, \emph{i.e.} $\tilde{\varphi}|_M = \varphi$.

Let $T\tilde{\varphi}$ be the tangent map of $\tilde{\varphi}$. It is a section of the bundle $\Omega^1(\overline{X})\otimes\tilde{\varphi}^* TN$, and its norm is defined by $$\|T\tilde{\varphi}\|^2_{g,h}:=\mathrm{tr}_g (\tilde{\varphi}^* h).$$ The bundle $\Omega^1(\overline{X})\otimes\tilde{\varphi}^* TN$ is canonically equipped with the connection $$\nabla^{g,h}:=\nabla^g\otimes 1_{\tilde{\varphi}^* TN}+1_{\Omega^1(\overline{X})}\otimes \nabla^{\tilde{\varphi}^*h},$$ where $\nabla^g$ and $\nabla^h$ are the respective Levi-Civita connections of $g$ and $h$, and $\nabla^{\tilde{\varphi}^*h}:=\tilde{\varphi}^*\nabla^{h}$.

\noindent The \emph{divergence} $\delta^{g,h}$ is then defined for $\omega\in \Omega^1(\overline{X})\otimes\tilde{\varphi}^* TN$ by $$\delta^{g,h}\omega:=-\left(\nabla^{g,h}_{e_I}\omega\right)(e_{\overline{I}}),$$ where $(e_i)$ is an orthonormal basis of $T^{1,0}\overline{X}$ for $g$, considered as a Hermitian metric.
We thus have $$\delta^{g,h}\omega=-\nabla^{\tilde{\varphi}^*h}_{e_I}(\omega(e_{\overline{I}}))+\omega(\nabla^g_{e_I} e_{\overline{I}}).$$

For $\rho\in(0,\varepsilon)$, the \emph{energy} of $\tilde{\varphi}$ in $(\rho,\varepsilon)\times M$ is the functional $$E(\tilde{\varphi},\rho)=\frac{1}{2}\int_{(\rho,\varepsilon)\times M}\|T\tilde{\varphi}\|^2_{g,h} d\mathrm{vol}_g.$$

\noindent An extension $\tilde{\varphi}$ is said to be \emph{harmonic} if it is a critical point of the energy for all $\rho$. Equivalently, $\tilde{\varphi}$ is harmonic if and only if $\delta^{g,h}T \tilde{\varphi} = 0$.

Following the ideas of C. R. Graham, R. Jenne, L. J. Mason, and G. A. J. Sparling, we want to find the obstructions to the existence of a smooth harmonic extension \cite{GJMS92}. More precisely, assuming that $\tilde{\varphi}$ is smooth, we want to know if the first terms of the asymptotic development of $\tilde{\varphi}$ are determined by the data at infinity. By similarity with the real case and based on the known asymptotic developments of the approximately ACH-Einstein metrics, we expect to find an obstruction at order $n+1$, taking the form of a CR covariant differential operator of order $2n+2$.

Here, the \emph{asymptotic development} of $\tilde{\varphi}$ will denote, by identification, the asymptotic development in $r$ of $U:=\exp_\varphi^{-1}\circ\tilde{\varphi}\in C^\infty(\overline{X},(\varphi\circ\pi)^* TN)$, \emph{i.e.} $$\forall p\in M,\ \forall r\in (0,\varepsilon),\quad \tilde{\varphi}(p,r):=\exp_{\varphi(p)}\left(U(p,r)\right),$$ where, for $p\in M$, the exponential map $\exp_{\varphi(p)}$ is a diffeomorphism between a small ball $B(0,\varepsilon)\subset T_{\varphi(p)}N$ and its image, which is a neighbourhood in $N$ of $\varphi(p)$.
Note that $U(\cdot,0)=0$. We denote $v\tilde{\varphi}:=T\tilde{\varphi}(v)$ for $v\in TX$, and similarly for $\varphi$ on $TM$, and $$\forall k\geq 1,\quad\varphi_k:=(\nabla^{\tilde{\varphi}^* h}_{\partial_r})^{k-1}\partial_r\tilde{\varphi}|_{r=0}.$$ Note that $\varphi_k$ is a section of $\varphi^*TN$, hence $\nabla^{\varphi^* h}\varphi_k$ is a section of $\Omega^1(M)\otimes\varphi^*TN$.

%Similarly, the \emph{horizontal divergence} $\delta^{W_r,h}$ is defined by $$\delta^{W_r,h}\omega=-\nabla^h_{T\tilde{\varphi}(e_i)}(\omega(e_i))+\omega(\nabla^\theta_{e_i} e_i)$$ where $(e_i)$ is an orthonormal basis of $\lbrace r\rbrace\times H$ with respect to $\gamma_r:=rg|_{\lbrace r\rbrace\times H}$.

%& \Longleftrightarrow k(2n+2-k)\partial_r^k\tilde{\varphi} + k(k-1)\partial_r^{k-2}\left(\Delta_{g_r}\tilde{\varphi}-\frac{1}{2}\partial_r\left(\log \det g_r\right)\partial_r\tilde{\varphi}\right)=0\\
%& \Leftrightarrow (n+1-k)\left(\nabla^{\tilde{\varphi}^*h}_{\partial_r}\right)^{k-1}T\tilde{\varphi}(\partial_r)|_{r=0} = -(k-1)\partial_r^{k-2}\left.\left(\Delta^{(1)}_{g_r}\tilde{\varphi}-\frac{1}{2}\mathrm{tr}^{(1)}_{g_r}(g'_r)\partial_r\tilde{\varphi}\right)\right|_{r=0}.

\subsection{Computation of the divergence}

We use the notations of section \ref{ache}. Let $(T_\alpha)$ be a local basis of $T^{1,0}M$ and $T_{\overline{\alpha}}:=\overline{T_\alpha}$, such that $(T_A)$ is orthonormal for $\gamma$, considered as a Hermitian metric. Let $(\theta^A)$ be the basis dual to $(T_A)$. Let $T_0:=\frac{\partial_r-i R}{\sqrt{2}}$ and $\theta^0:=\frac{dr+i\theta}{\sqrt{2}}$ its dual.

\begin{lemm}\label{comput}
For $\omega\in \Omega^1(\overline{X})\otimes\tilde{\varphi}^* TN$, we have
\begin{displaymath}
\begin{array}{lll}
{\delta^{g_0,h}\omega}&=&nr\omega(\partial_r)-r^2\left(\nabla^{\tilde{\varphi}^*h}_{T_0}\omega(T_{\overline{0}})+\nabla^{\tilde{\varphi}^*h}_{T_{\overline{0}}}\omega(T_0)\right)-r\nabla^{\tilde{\varphi}^*h}_{T_A}\omega(T_{\overline{A}})\\
&=& nr\omega(\partial_r)-r^2\nabla^{\tilde{\varphi}^*h}_{\partial_r}\omega(\partial_r)-r^2\nabla^{\tilde{\varphi}^*h}_{R}\omega(R)-r\nabla^{\tilde{\varphi}^*h}_{T_A}\omega(T_{\overline{A}}).
\end{array}
\end{displaymath}
\end{lemm}

\begin{proof}

We have $$g_0=r^{-2}\theta^0\circ\theta^{\overline{0}} + r^{-1}\theta^\alpha\circ\theta^{\overline{\alpha}}.$$%  where $\lambda\circ\mu:=\lambda\otimes\mu+\mu\otimes\lambda$.

An orthonormal basis of $T^{1,0}\overline{X}$ with respect to $g_0$ is hence given by
\begin{displaymath}
\begin{array}{lll}
(e_0^{(0)},e_\alpha^{(0)})&:=&\left(rT_0, r^\frac{1}{2}T_\alpha\right).
\end{array}
\end{displaymath}
%Let $O_e$ denote the order with respect to $(e_i^{(0)})$. 

%Now, $g$ can be rewritten as $$g=(r^{-1}\theta^0)\circ(r^{-1}\theta^{\overline{0}})+(r^{-\frac{1}{2}}\theta^\alpha)\circ(r^{-\frac{1}{2}}\theta^{\overline{\alpha}})+r\Phi_{AB}(r^{-\frac{1}{2}}\theta^A)\circ(r^{-\frac{1}{2}}\theta^B)+O_e(r^\frac{3}{2}).$$

\noindent The trace of the Levi-Civita connection of $g_0$ is given in this basis by the Koszul formula: $$\nabla^{g_0}_{e_I^{(0)}}e_{\overline{I}}^{(0)}=g_0\left([e_{\overline{J}}^{(0)},e_I^{(0)}],e_{\overline{I}}^{(0)}\right)e_J^{(0)}.$$
%Let us denote by $\rho$ the matrix of $T^W(R,\cdot)$ in the basis $h$: $$T^W(R,h_\alpha)=\rho_{\alpha\beta}h_\beta+\rho_{\alpha\tilde{\beta}}Jh_\beta,$$ $$T^W(R,Jh_\alpha)=\rho_{\tilde{\alpha}\beta}h_\beta+\rho_{\tilde{\alpha}\tilde{\beta}}Jh_\beta.$$ Since $T^W(R,\cdot)$ and $J$ anticommute, $\rho$ is block-symmetric and -trace-free, \emph{i.e.} $$\rho_{\tilde{\alpha}\beta}=\rho_{\alpha\tilde{\beta}}\quad\mathrm{and}\quad\rho_{\tilde{\alpha}\tilde{\beta}}=-\rho_{\alpha\beta}.$$
Let $\tilde{\nabla}^\theta$ be the extension of $\nabla^\theta$ given by Proposition \ref{comp2}. We have
\begin{displaymath}
\begin{array}{rcl}
[e_0^{(0)},e_{\overline{0}}^{(0)}]&=&\frac{1}{\sqrt{2}}\left(e_{\overline{0}}^{(0)}-e_0^{(0)}\right),\\

[e_0^{(0)},e_A^{(0)}]&=&\frac{1}{\sqrt{2}}\left(\frac{1}{2} e_A^{(0)}-i\left(\tilde{\nabla}^\theta_{{e_0}^{(0)}}e_A^{(0)}-\tau(e_A^{(0)})\right)\right),\\

%[e_0^{(0)},e_{\overline{\alpha}}^{(0)}]&=&\frac{1}{\sqrt{2}}\left(\frac{1}{2} e_{\overline{\alpha}}^{(0)}-i\left(\tilde{\nabla}^\theta_{{e_0}^{(0)}}e_{\overline{\alpha}}^{(0)}-\tau(e_{\overline{\alpha}}^{(0)})\right)\right),\\

[e_A^{(0)},e_B^{(0)}]&=&rd\theta(T_A,T_B)R.

\end{array}
\end{displaymath}
Then, since $\mathrm{tr}(\tau)=0$, $$\nabla^{g_0}_{e_I^{(0)}}e_{\overline{I}}^{(0)}=\left(n+1\right)r\partial_r,$$ and also,
\begin{displaymath}
\begin{array}{rcl}
\nabla^{\tilde{\varphi}^*h}_{e_0^{(0)}}\omega(e_{\overline{0}}^{(0)})+\nabla^{\tilde{\varphi}^*h}_{e_{\overline{0}}^{(0)}}\omega(e_0^{(0)})&=&r\omega(\partial_r)+r^2\nabla^{\tilde{\varphi}^*h}_{\partial_r}\omega(\partial_r)+r^2\nabla^{\tilde{\varphi}^*h}_R\omega(R),\\

\nabla^{\tilde{\varphi}^*h}_{e_\alpha^{(0)}}\omega(e_{\overline{\alpha}}^{(0)})&=&r\nabla^{\tilde{\varphi}^*h}_{T_\alpha}\omega(T_{\overline{\alpha}}).
\end{array}
\end{displaymath}

\noindent Hence the announced expression for $\delta^{g_0,h}\omega$.
\end{proof}

Let us denote by $\left({\delta^{g,h}\omega}\right)^{(1)}$ the remainder of $\delta^{g,h}\omega$, \emph{i.e.} $$\left({\delta^{g,h}\omega}\right)^{(1)}:=\delta^{g,h}\omega-{\delta^{g_0,h}\omega}.$$

We prove the following technical lemma, which is crucial for the proof of Theorem \ref{obs}.

\begin{lemm}\label{comput2}
For $\omega\in \Omega^1(\overline{X})\otimes\tilde{\varphi}^* TN$, denoting by $O_T$ the order with respect to the basis $(\tilde{\varphi}_* T_I)$ in powers of $r$, we have $$\left({\delta^{g,h}\omega}\right)^{(1)}=O_T(r^2),$$ and there is no term of order $2$ in the remainder of the form $r^2\nabla^{\tilde{\varphi}^*h}_{\partial_r}\omega(\partial_r)$.
\end{lemm}

\begin{proof}
By Theorem \ref{poinc}, we have $$g-g_0=\Phi+O_e(r^\frac{3}{2})=\Phi_{AB}\theta^A\circ\theta^B+O_e(r^\frac{3}{2}),$$ where we recall that $\Phi = -2 \mathrm{Sch}_W(J,\theta) + 2 \gamma(J\tau\cdot,\cdot)$, and that $O_e$ denotes the order with respect to $(e_I^{(0)})$. Note that $\Phi_{AB}=\Phi_{BA}$. Since $\Phi$ is real, we have also $\Phi_{\overline{\alpha}\beta}=\overline{\Phi_{\alpha\overline{\beta}}}$ and $\Phi_{\overline{\alpha}\overline{\beta}}=\overline{\Phi_{\alpha\beta}}$.

An orthonormal basis of $T^{1,0}\overline{X}$ with respect to $g$ induced from $e^{(0)}$ is formally given by
\begin{displaymath}
\begin{array}{lll}
(e_0,e_\alpha)&:=&\left(e_0^{(0)}+e_0^{(1)}, e_\alpha^{(0)}+e_\alpha^{(1)}\right),
\end{array}
\end{displaymath}
where, by the Gram-Schmidt process, and since $\Phi$ is horizontal, $$e_0^{(1)}=O_e(r^\frac{3}{2})\quad \text{and}\quad e_\alpha^{(1)}=O_e(r).$$
%The trace of the Levi-Civita connection of $g$ is, as in the previous lemma, given in the basis $(e_I)$ by the Koszul formula.

\noindent This leads to
\begin{displaymath}
\begin{array}{rcl}
\left(\delta^{g,h}\omega\right)^{(1)} &=& 
-\nabla^{\tilde{\varphi}^*h}_{e_I^{(0)}}\omega\left(e_{\overline{I}}^{(1)}\right)-\nabla^{\tilde{\varphi}^*h}_{e_I^{(1)}}\omega\left(e_{\overline{I}}^{(0)}\right)-\nabla^{\tilde{\varphi}^*h}_{e_I^{(1)}}\omega\left(e_{\overline{I}}^{(1)}\right)\\
&&+\omega\left(\nabla^g_{e_I^{(0)}}e_{\overline{I}}^{(0)}-\nabla^{g_0}_{e_I^{(0)}}e_{\overline{I}}^{(0)}\right)+ \omega\left(\nabla^{g}_{e_I^{(0)}}e_{\overline{I}}^{(1)}\right)+ \omega\left(\nabla^{g}_{e_I^{(1)}}e_{\overline{I}}^{(0)}\right)+ \omega\left(\nabla^{g}_{e_I^{(1)}}e_{\overline{I}}^{(1)}\right),
\end{array}
\end{displaymath}
all terms of which are in $O_T(r^2)$ and are not of the form $r^2\nabla^{\tilde{\varphi}^*h}_{\partial_r}\omega(\partial_r)$.
\end{proof}

\subsection{An obstruction to regularity}

\begin{theo}\label{obs}
Let $(M,H,J)$ be a $(2n+1)$-dimensional orientable, compact, strictly pseudoconvex CR manifold and $(X,g)$ be an ACH manifold with CR infinity $(M,H,J)$, where $g$ is the approximately ACH-Kähler-Einstein metric given by Theorem \ref{ke}. Let $\pi:\overline{X}\rightarrow M$ be the standard projection. Let $(N,h)$ be a Riemannian manifold, and let $\varphi \in C^\infty(M,N)$.

There exists a section $U$ of $(\varphi\circ\pi)^* TN$, unique modulo $O_T(r^{n+1})$, such that $\tilde{\varphi}=\exp_\varphi\circ U$ satisfies
\begin{displaymath}
\left\lbrace\begin{array}{rcl}
\tilde{\varphi}|_M&=&\varphi,\\
\delta^{g,h}T\tilde{\varphi}&=&O_T(r^{n+2}).
\end{array}\right.
\end{displaymath}

The asymptotic development in $r$ of $U$ is
$$U = U_1 r +\ldots+ U_{n}\frac{r^n}{n!}+P_n(\varphi)\frac{r^{n+1}}{(n+1)!}\log r +O_T(r^{n+1}),$$
where $U_1,\ldots,U_n,P_n$ are formally determined by $\varphi$, $g$ and $h$.

$P_n(\varphi)$ is an obstruction to the regularity of $U$, and is given by
\begin{displaymath}
\begin{array}{rcl}

P_n(\varphi)&=& \left.\left(\nabla^{\tilde{\varphi}^*h}_{\partial_r}\right)^n\tilde{\delta}^{\theta,h}_b T\tilde{\varphi}\right|_{r=0} -n\left.\left(\nabla^{\tilde{\varphi}^*h}_{\partial_r}\right)^{n-1} \nabla^{\tilde{\varphi}^*h}_R R\tilde{\varphi}\right|_{r=0}+\frac{1}{n+1}\left.\left(\nabla^{\tilde{\varphi}^*h}_{\partial_r}\right)^{n+1}{\left(\delta^{g,h}T\tilde{\varphi}\right)}^{(1)}\right|_{r=0}\\
&=&\frac{(-1)^n}{n!}(\delta^{\theta,h}_b\nabla^{\varphi^*h})^n \delta^{\theta,h}_bT\varphi+\text{lower order terms (in derivatives of $\varphi$)}.

\end{array}
\end{displaymath}
\end{theo}

\begin{proof}

For $m\in\N$, we have
$$\delta^{g,h}T\tilde{\varphi}=O_T(r^{m+1}) \quad\Longleftrightarrow\quad \forall k\leq m,\quad \left.\left(\nabla^{\tilde{\varphi}^*h}_{\partial_r}\right)^k\delta^{g,h}T\tilde{\varphi}\right|_{r=0}=0.$$
We recall the notation $$\varphi_k:=\left.\left(\nabla^{\tilde{\varphi}^* h}_{\partial_r}\right)^{k-1}\partial_r\tilde{\varphi}\right|_{r=0}.$$ Now, by Lemma \ref{comput}, we have, for $\omega\in\Omega^1(\overline{X})\otimes\tilde{\varphi}^*TN$,
$$\left.\nabla^{\tilde{\varphi}^*h}_{\partial_r} \delta^{g,h}\omega\right|_{r=0}=n\left.\omega(\partial_r)\right|_{r=0}+\delta^{\theta,h}_b (\omega|_{r=0}),$$ and, for all $2\leq k\leq n$,
\begin{displaymath}
\begin{array}{rcl}
\frac{1}{k}\left.\left(\nabla^{\tilde{\varphi}^*h}_{\partial_r}\right)^k \delta^{g,h}\omega\right|_{r=0}&=& (n-k+1)\left.\left(\nabla^{\tilde{\varphi}^*h}_{\partial_r}\right)^{k-1}\omega(\partial_r)\right|_{r=0}+\left.\left(\nabla^{\tilde{\varphi}^*h}_{\partial_r}\right)^{k-1}
\tilde{\delta}^{\theta,h}_b \omega\right|_{r=0}\\
&&-(k-1)\left.\left(\nabla^{\tilde{\varphi}^*h}_{\partial_r}\right)^{k-2}\nabla^{\tilde{\varphi}^*h}_R \omega(R)\right|_{r=0}+\frac{1}{k}\left.\left(\nabla^{\tilde{\varphi}^*h}_{\partial_r}\right)^k{\left(\delta^{g,h}\omega\right)}^{(1)}\right|_{r=0},
\end{array}
\end{displaymath}
where $$\forall \omega_0\in \Omega^1(M)\otimes\varphi^*TN,\quad \delta^{\theta,h}_b \omega_0:=-\nabla^{\varphi^*h}_{T_A}\omega_0(T_{\overline{A}}),$$ and $$\forall \omega\in \Omega^1(\overline{X})\otimes\tilde{\varphi}^*TN,\quad \tilde{\delta}^{\theta,h}_b \omega:=-\nabla^{\tilde{\varphi}^*h}_{T_A}\omega(T_{\overline{A}}).$$

Consequently, $\delta^{g,h}T\tilde{\varphi}=O_T(r^{n+1})$ is equivalent to
\begin{displaymath}
\left\lbrace
\begin{array}{rcl}
n\varphi_1&=&-\delta^{\theta,h}_b T\varphi,\\
(n-k+1)\varphi_k&=&-\left.\left(\nabla^{\tilde{\varphi}^*h}_{\partial_r}\right)^{k-1}\tilde{\delta}^{\theta,h}_b T\tilde{\varphi}\right|_{r=0}-D_{k-1}(\varphi)\quad\forall 2\leq k \leq n,
\end{array}
\right.
\end{displaymath}
where $$D_{k-1}(\varphi):=-(k-1)\left.\left(\nabla^{\tilde{\varphi}^*h}_{\partial_r}\right)^{k-2} \nabla^{\tilde{\varphi}^*h}_R R\tilde{\varphi}\right|_{r=0}+\frac{1}{k}\left.\left(\nabla^{\tilde{\varphi}^*h}_{\partial_r}\right)^k{\left(\delta^{g,h}T\tilde{\varphi}\right)}^{(1)}\right|_{r=0}.$$
By Lemma \ref{comput2}, $D_{k-1}(\varphi)$ only depends on $\varphi,\varphi_1,\ldots,\varphi_{k-1}$. This observation comes from the fact that, although $$\left.\left(\nabla^{\tilde{\varphi}^*h}_{\partial_r}\right)^k \left(r^2\nabla^{\tilde{\varphi}^*h}_{\partial_r}\partial_r\tilde{\varphi}\right)\right|_{r=0}=2\varphi_k,$$ 
$$\forall X,Y\in\lbrace \partial_r,R,T_A\rbrace,\ (X,Y)\neq(\partial_r,\partial_r),\quad \left.\left(\nabla^{\tilde{\varphi}^*h}_{\partial_r}\right)^k \left(r^2\nabla^{\tilde{\varphi}^*h}_X Y\tilde{\varphi}\right)\right|_{r=0} \ \text{
does not depend on $\varphi_k$.}$$ By induction, $D_{k-1}(\varphi)$ is thus well-defined.

In conclusion, requiring $\delta^{g,h}T\tilde{\varphi}=O_T(r^{n+1})$ gives an asymptotic development for $\tilde{\varphi}$ in powers of $r$, and this development is unique up to order $n$ with respect to $T$.
%and 
%\begin{displaymath}
%\begin{array}{ll}
%\square \omega:=& -\nabla^{\tilde{\varphi}^*h}_R(\omega(R))+2\mathrm{Im}\left(\nabla^{\tilde{\varphi}^*h}_{T_{\overline{\alpha}}}(A_\alpha^{\overline{\beta}}+A_\beta^{\overline{\alpha}})\omega(T_{\overline{\beta}})\right)\\
%&+\frac{1}{2}\nabla^{\tilde{\varphi}^*h}_{h_\alpha}(\Phi_{\alpha\overline{\beta}}+\Phi_{\beta\overline{\alpha}})\omega(h_\beta)-\nabla^{\tilde{\varphi}^*h}_{h_\alpha}\Phi_{\beta\overline{\beta}}\omega(h_\alpha)+C\omega(h_\alpha)\\
%&+\frac{1}{2}\nabla^{\tilde{\varphi}^*h}_{J_0h_\alpha}(\Phi_{\alpha\overline{\beta}}+\Phi_{\beta\overline{\alpha}})\omega(J_0h_\beta)-\nabla^{\tilde{\varphi}^*h}_{J_0h_\alpha}\Phi_{\beta\overline{\beta}}\omega(J_0h_\alpha)+D\omega(J_0 h_\alpha).
%\end{array}
%\end{displaymath}

Assume now that $\delta^{g,h}T\tilde{\varphi}=O_T(r^{n+1})$ and that $\tilde{\varphi}$ admits a Taylor development up to order $n+1$. Then %there is no term of order $2n+1$ in the development of $X$ $\partial_r^{n+1}\tilde{\varphi}|_{r=0}$ exists, and 
$$\delta^{g,h}T\tilde{\varphi}=O_T(r^{n+2})\Longleftrightarrow \left.\left(\nabla^{\tilde{\varphi}^*h}_{\partial_r}\right)^n\tilde{\delta}^{\theta,h}_b T\tilde{\varphi}\right|_{r=0}+D_n(\varphi)=0.$$ This equality cannot be true in general. Consequently, we introduce a term in $r^{n+1}\log r$ in the development of $\tilde{\varphi}$:
$$U = U_1 r +\ldots+ U_{n}\frac{r^n}{n!}+P_n(\varphi)\frac{r^{n+1}}{(n+1)!}\log r +O_T(r^{n+1}).$$
The coefficient $P_n(\varphi)$ verifies $$\frac{1}{n+1}\left.\left(\nabla^{\tilde{\varphi}^*h}_{\partial_r}\right)^{n+1} \delta^{g,h}T\tilde{\varphi}\right|_{r=0}= -P_n(\varphi)+\left.\left(\nabla^{\tilde{\varphi}^*h}_{\partial_r}\right)^n\tilde{\delta}^{\theta,h}_b T\tilde{\varphi}\right|_{r=0}+D_n(\varphi),$$
hence $$\delta^{g,h}T\tilde{\varphi}=O_T(r^{n+2})\Longleftrightarrow P_n(\varphi)=\left.\left(\nabla^{\tilde{\varphi}^*h}_{\partial_r}\right)^n\tilde{\delta}^{\theta,h}_b T\tilde{\varphi}\right|_{r=0}+D_n(\varphi).$$

\noindent This yields the announced obstruction, which only depends on $\varphi$. Since $$\varphi_k=-\frac{1}{n-k+1}\delta^{\theta,h}_b\nabla^{\varphi^* h}\varphi_{k-1}+\text{lower order terms (in derivatives of $\varphi$)},$$ we have the announced leading term.
\end{proof}

\begin{prop}\label{suff}
$P_n$ does not depend on whether we take $g=g_E$ or $g_{KE}$
%we can equivalently consider the approximately ACH-Einstein metric or the approximately Kähler ACH-Einstein 
 on $X$.
\end{prop}

\begin{proof}
To compute $P_n$, it is sufficient to be able to compute $$\left.\left(\nabla^{\tilde{\varphi}^*h}_{\partial_r}\right)^{n+1}\left(\delta^{g,h}T\tilde{\varphi}\right)^{(1)}\right|_{r=0};$$ \emph{i.e.}, by the proof of Lemma \ref{comput2}, to know the $e_I^{(1)}$ at order $n+1/2$ with respect to $e^{(0)}$. By the Gram-Schmidt process, it is thus sufficient to know $g$ at order $n+1/2$ with respect to $e^{(0)}$. Hence, by Theorems \ref{poinc} and \ref{ke}, we can equivalently consider $g_E$ or $g_{KE}$.
\end{proof}

%As an example, let us consider $\mathrm{id}:M\rightarrow M$. We have $$P_n(\mathrm{id})=\frac{1}{n+1}\left.\left(\nabla^{\tilde{\varphi}^*h}_{\partial_r}\right)^{n+1} {\left(\delta^{g,h}T\tilde{\varphi}\right)}^{(1)}\right|_{r=0}.$$

\begin{prop}\label{crinv}
Let $f\in C^\infty(\overline{X},\R)$ and $f_0:=f\vert_M$, and let $\hat{r}=e^{f} r$ be a conformal change of boundary defining function. Then $$\hat{P}_n(\varphi)=e^{-(n+1)f_0}P_n(\varphi).$$
\end{prop}
\noindent The obstruction $P_n(\varphi)$ to the regularity of $\tilde{\varphi}$ is therefore CR covariant.

\begin{proof}
We have $$ U= U_1 r +\ldots+ U_{n}\frac{r^n}{n!}+P_n(\varphi)\frac{r^{n+1}}{(n+1)!}\log r +O_T(r^{n+1}).$$
%Now, let $\hat{r}=e^f r$ be a conformal change of boundary defining function. Then,
%Note that $U$ does not depend on the choice of $r$. In particular, 
Now, since $\exp_\varphi:(\varphi\circ\pi)^*TN\rightarrow N$ does not depend on $r$, neither does $U$. Moreover, since $M$ is compact, %$\forall r,\hat{r},\ 
$\forall k,\ O_T(\hat{r}^k)=O_T(r^k)$. We thus have
\begin{displaymath}
\begin{array}{rcl}
U & = & \hat{U}_1 \hat{r} +\ldots+ \hat{U}_{n}\frac{\hat{r}^{n}}{n!}+\hat{P}_n(\varphi)\frac{\hat{r}^{n+1}}{(n+1)!}\log \hat{r} +O_T(\hat{r}^{n+1})\\
& = & \hat{U}_1 e^f r +\ldots+ \hat{U}_{n}e^{nf} \frac{r^n}{n!}+\hat{P}_n(\varphi)e^{(n+1)f}\frac{r^{n+1}}{(n+1)!}\log r +O_T(r^{n+1}).
\end{array}
\end{displaymath}
Since the function $f$ itself has a Taylor expansion in $r$, all polynomial terms are mixed. However, there is only one term with order $r^{n+1}\log r$. By identification, this yields the result.
\end{proof}

We then introduce \emph{CR-harmonic maps} as maps for which the obstruction vanishes:

\begin{defi}
If $P_n(\varphi)=0$, $\varphi$ is said to be \emph{CR-harmonic}.
\end{defi}

%\emph{Remark:} If $\nabla^{\tilde{\varphi}^*h}_R R\tilde{\varphi} =-\nabla^{\tilde{\varphi}^*h}_{\partial_r} \partial_r\tilde{\varphi}$, then $P_n(\varphi)=0$.\\ %($N$ complex ? $J_0\partial_r=R$ by definition)

\begin{ex}\label{psn}
Let us assume that $(M,H,J,\theta)$ is Einstein with $\mathrm{Ric}_W=2\lambda(n+1)\gamma$. %Distinction 3/others : cf Case-Yang p.11
We know from Example \ref{exlam} that $$g = \frac{dr^2}{r^2}+\frac{(1-\lambda^2r^2)^2}{r^2}\theta^2+\frac{(1-\lambda r)^2}{r}\gamma$$ satisfies $\mathrm{Ric}(g)=-\frac{n+2}{2}g$. %Consequently, if $n$ is even, subharmonic maps are CR-harmonic.
In this case, we can explicitly compute the divergence $\delta^{g,h}\omega$, for $\omega\in\Omega^1(\overline{X})\otimes\tilde{\varphi}^*TN$.
%\begin{proof}

\noindent Indeed, an orthonormal basis of $T^{1,0}\overline{X}$ with respect to $g$ induced from $e^{(0)}$ is given by
\begin{displaymath}
\begin{array}{lll}
(e_0,e_\alpha)&:=&\left(\frac{1}{\sqrt{2}}\left(r\partial_r-i\frac{r}{1-\lambda^2r^2}R\right), \frac{r^\frac{1}{2}}{1-\lambda r}T_\alpha\right),
\end{array}
\end{displaymath}
hence
\begin{displaymath}
\begin{array}{rcl}
[e_0,e_{\overline{0}}]&=&\frac{1}{\sqrt{2}}\frac{1+\lambda^2r^2}{1-\lambda^2r^2}\left(e_{\overline{0}}-e_0\right),\\

[e_0,e_A]&=&\frac{1}{2\sqrt{2}}\frac{1+\lambda r}{1-\lambda r} e_A,\\

%[e_0,e_{\overline{\alpha}}]&=&\frac{1}{2\sqrt{2}}\frac{1+\lambda r}{1-\lambda r} e_{\overline{\alpha}},\\

[e_A,e_B]%=[e_\alpha,e_{\overline{\beta}}]
&=&\frac{r}{(1-\lambda r)^2}d\theta(T_A,T_B)R.
\end{array}
\end{displaymath}

\noindent Then $$\nabla^{g}_{e_I}e_{\overline{I}}=\left(n\frac{1+\lambda^2r^2}{1-\lambda^2r^2}+\frac{1+\lambda r}{1-\lambda r}\right)r\partial_r,$$ and also,
\begin{displaymath}
\begin{array}{rcl}
\nabla^{\tilde{\varphi}^*h}_{e_0}\omega(e_{\overline{0}})+\nabla^{\tilde{\varphi}^*h}_{e_{\overline{0}}}\omega(e_0)&=&r\omega(\partial_r)+r^2\nabla^{\tilde{\varphi}^*h}_{\partial_r}\omega(\partial_r)+\frac{r^2}{(1-\lambda^2r^2)^2}\nabla^{\tilde{\varphi}^*h}_R\omega(R),\\

\nabla^{\tilde{\varphi}^*h}_{e_\alpha}\omega(e_{\overline{\alpha}})&=&\frac{r}{(1-\lambda r)^2}\nabla^{\tilde{\varphi}^*h}_{T_\alpha}\omega(T_{\overline{\alpha}}).
\end{array}
\end{displaymath}

\noindent The divergence is hence given by
\begin{displaymath}
\begin{array}{lll}
\delta^{g,h}\omega&=&\left(n\frac{1+\lambda^2r^2}{1-\lambda^2r^2}+\frac{1+\lambda r}{1-\lambda r}-1\right)r\omega(\partial_r)-r^2\nabla^{\tilde{\varphi}^*h}_{\partial_r}\omega(\partial_r)\\
&&-\frac{r^2}{(1-\lambda^2r^2)^2}\nabla^{\tilde{\varphi}^*h}_R\omega(R)+\frac{r}{(1-\lambda r)^2}\tilde{\delta}^{\theta,h}_b\omega.
\end{array}
\end{displaymath}
%Let us assume that $n$ is even and that $\varphi_1=-\frac{1}{n}\delta^{\theta,h}_bT\varphi=0$. Then, for all $2\leq k\leq n$, $\varphi_2=\ldots=\varphi_{k-1}=0$ implies $$\frac{1}{k}\left.\left(\nabla^{\tilde{\varphi}^*h}_{\partial_r}\right)^k \delta^{g,h}T\tilde{\varphi}\right|_{r=0}= (n-k+1)\varphi_k,$$ hence, by induction, $P_n(\varphi)=0$.
%\end{proof}

\end{ex}

From Example \ref{psn} we get the following results:

\begin{cor}\label{sub}
If $(M,H,J,\theta)$ is Einstein, then subharmonic maps which verify $\nabla^{\varphi^*h}_R R\varphi=0$ are CR-harmonic.
\end{cor}
\begin{proof}
Indeed, let $\varphi$ be subharmonic, \emph{i.e.} $\delta^{\theta,h}_bT\varphi=0$, and such that $\nabla^{\varphi^*h}_R R\varphi=0$. Let $\tilde{\varphi}$ be the extension of $\varphi$ given by Theorem \ref{obs}. We thus have $\varphi_1=0$. Moreover, by Example \ref{psn}, we have $${\left(\delta^{g,h}T\tilde{\varphi}\right)}^{(1)}=\alpha(r)\partial_r\tilde{\varphi}+\beta(r)\nabla^{\tilde{\varphi}^* h}_R R\tilde{\varphi}+\gamma(r)\tilde{\delta}^{\theta,h}_bT\tilde{\varphi},$$ where $\alpha(r)=O(r^2)$, $\beta(r)=O(r^4)$, and $\gamma(r)=O(r^2)$. Since $\varphi_1=\nabla^{\varphi^*h}_R R\varphi=0$, we get that $$(n-1)\varphi_2=-\left.\nabla^{\tilde{\varphi}^* h}_{\partial_r}\tilde{\delta}^{\theta,h}_bT\tilde{\varphi}\right|_{r=0}-\nabla^{\varphi^*h}_R R\varphi-\frac{1}{2}\left.\left(\nabla^{\tilde{\varphi}^* h}_{\partial_r}\right)^2{\left(\delta^{g,h}T\tilde{\varphi}\right)}^{(1)}\right|_{r=0} =0. $$ By induction, we similarly have $\forall k\leq n,\ \varphi_k=0$, which implies that $P_n(\varphi)=0$.%$$\delta^{g,h}T\tilde{\varphi}=O_T(r^{k+1})\Longleftrightarrow \forall j\leq k,\ \varphi_j=0,$$ 
%\begin{displaymath}
%\begin{array}{rcl}
%0&=&\int_M g_{J,\theta}\left(\delta^{\theta,h}_bT\varphi,\overline{\varphi}\right)\theta\wedge d\theta^n\\
%&=&\int_M \left(g_{J,\theta}\left(T_\alpha\varphi,T_{\overline{\alpha}}\overline{\varphi}\right)+g_{J,\theta}\left(T_{\overline{\alpha}}\varphi,T_\alpha\overline{\varphi}\right)\right)\theta\wedge d\theta^n,
%\end{array}
%\end{displaymath}
%hence, $\forall 1\leq \alpha\leq n$, $T_\alpha\varphi=T_{\overline{\alpha}}\varphi=0$, hence $R\varphi=0$ by the following fact:
%\begin{fct}\label{brackets}
%We have% (cf GHL 1.66)
% $$\nabla^{\varphi^*h}_{T_\alpha} T_{\overline{\alpha}}\varphi-\nabla^{\varphi^*h}_{T_{\overline{\alpha}}} T_\alpha\varphi=\varphi_*[T_\alpha,T_{\overline{\alpha}}]=iR\varphi,$$ and $$\nabla^{\varphi^*h}_R T_{\overline{\alpha}}\varphi-\nabla^{\varphi^*h}_{T_{\overline{\alpha}}} R\varphi=\varphi_*[R,T_{\overline{\alpha}}]=\tau_{\overline{\alpha}}^\beta T_\beta\varphi.$$
%\end{fct}
% The divergence is then given from Example \ref{psn} by $$\delta^{g,h}\omega=\left(n\frac{1+\lambda^2r^2}{1-\lambda^2r^2}+\frac{1+\lambda r}{1-\lambda r}-1\right)r\omega(\partial_r)-r^2\nabla^{\tilde{\varphi}^*h}_{\partial_r}\omega(\partial_r)+\frac{r}{(1-\lambda r)^2}\tilde{\delta}^{\theta,h}_b\omega.$$ From this we easily see by induction that $$\delta^{g,h}T\tilde{\varphi}=O_T(r^{k+1})\Longleftrightarrow \forall j\leq k,\ \varphi_j=0,$$ hence $P_1(\varphi)=0$.
\end{proof}

\begin{cor}
If $(M,H,J,\theta)$ is Einstein and $(N,h)$ is a Kähler manifold, then CR-holomorphic maps which verify $R\varphi=0$ are CR-harmonic.
\end{cor}
\begin{proof}
Indeed, %let us denote $h_\alpha:=\frac{T_\alpha+T_{\overline{\alpha}}}{\sqrt{2}}$ and $Jh_\alpha:=J(h_\alpha)$. We 
assuming that $T\varphi\circ J=J_N\circ T\varphi$, and extending $J$ by taking $J(R)=0$, we have
\begin{displaymath}
\begin{array}{rcl}
\nabla^{\varphi^*h}_{T_\alpha}T_{\overline{\alpha}}\varphi&=&\nabla^{\varphi^*h}_{JT_\alpha}JT_{\overline{\alpha}}\varphi\\
&=&J_N\nabla^{\varphi^*h}_{JT_\alpha}T_{\overline{\alpha}}\varphi\\
&=&J_N\nabla^{\varphi^*h}_{T_{\overline{\alpha}}}JT_\alpha\varphi+J([JT_\alpha,T_{\overline{\alpha}}])\varphi\\
&=&-\nabla^{\varphi^*h}_{T_{\overline{\alpha}}}T_\alpha\varphi+iJ([T_\alpha,T_{\overline{\alpha}}])\varphi\\
&=&-\nabla^{\varphi^*h}_{T_{\overline{\alpha}}}T_\alpha\varphi-nJ(R)\varphi,
%&=&-\nabla^{\varphi^*h}_{T_{\overline{\alpha}}}T_\alpha\varphi,
\end{array}
\end{displaymath}
%\begin{displaymath}
%\begin{array}{rcl}
%\nabla^{\varphi^*h}_{Jh_\alpha}Jh_\alpha\varphi&=&J_N\nabla^{\varphi^*h}_{Jh_\alpha}h_\alpha\varphi\\
%&=&J_N\nabla^{\varphi^*h}_{h_\alpha}Jh_\alpha\varphi-J([h_\alpha,Jh_\alpha])\varphi\\
%&=&-\nabla^{\varphi^*h}_{h_\alpha}h_\alpha\varphi+iJ([T_\alpha,T_{\overline{\alpha}}])\varphi\\
%&=&-\nabla^{\varphi^*h}_{h_\alpha}h_\alpha\varphi-J(R)\varphi\\
%&=&-\nabla^{\varphi^*h}_{h_\alpha}h_\alpha\varphi,
%\end{array}
%\end{displaymath}
hence $$\delta^{\theta,h}_bT\varphi=nJ_N\left(R\varphi\right).$$ Consequently, $\varphi$ is CR-harmonic by Corollary \ref{sub}.
\end{proof}

\begin{ex}
Let $(M,H,J)$ be a circle bundle over a Riemann surface $\Sigma$ admitting an Einstein contact form. Then the projection $\pi:M\rightarrow\Sigma$ is CR-harmonic.
\end{ex}

\subsection{Explicit obstruction in dimension $3$}
\label{exob}

When $n=1$, \emph{i.e.} $\dim(M)=3$, the asymptotic development of $g$ is given at order $\frac{3}{2}$ in $e^{(0)}$ by Theorem \ref{expl}. Hence, by Proposition \ref{suff}, we can explicitly compute the obstruction.

\begin{theo}\label{corpsB1}
Still denoting $v\varphi:=T\varphi(v)$ for $v\in TM$, and also $(\nabla^{\varphi^* h}v)\varphi:=\nabla^{\varphi^* h}(v\varphi)$, we have $$P_1(\varphi)=-\delta^{\theta,h}_b \nabla^{\varphi^*h}\delta^{\theta,h}_bT\varphi-\nabla^{\varphi^* h}_R R\varphi+4\mathrm{Im}\left(\nabla^{\varphi^* h}_{T_{\overline{1}}}\left(\tau_1^{\overline{1}} T_{\overline{1}}\right)\right)\varphi-S_b\left(\delta^{\theta,h}_bT\varphi\right),$$
where $$S_b(X):=\mathcal{R}^h_{X,T_1\varphi}T_{\overline{1}}\varphi+\mathcal{R}^h_{X,T_{\overline{1}}\varphi}T_1\varphi.$$
%$$\quad S(X):=R^h_{X,T_\alpha\varphi}T_{\overline{\alpha}}\varphi+R^h_{X,T_{\overline{\alpha}}\varphi}T_\alpha\varphi.$$
\end{theo}

\begin{proof}
 
An orthonormal basis of $T^{1,0}X$ with respect to $g$ is given by
\begin{displaymath}
\begin{array}{ll}
(e_0,e_1):=&\left(e_0^{(0)}-r^\frac{3}{2}\Psi_{0\overline{1}}e_1^{(0)},\left(1 -r\Phi_{1\overline{1}}\right)e_1^{(0)} -r\Phi_{11} e_{\overline{1}}^{(0)}\right)+O_e(r^2).
\end{array}
\end{displaymath}
%We can thus rewrite
%\begin{displaymath}
%\begin{array}{rcl}
%g&=&(r^{-1}\theta^0)\circ(r^{-1}\theta^{\overline{0}})+(r^{-\frac{1}{2}}\theta^1)\circ(r^{-\frac{1}{2}}\theta^{\overline{1}})\\\\
%&&+r\left(2\Phi_{1\overline{1}}(r^{-\frac{1}{2}}\theta^1)\circ(r^{-\frac{1}{2}}\theta^{\overline{1}})+\Phi_{11}(r^{-\frac{1}{2}}\theta^1)\circ(r^{-\frac{1}{2}}\theta^1)+\Phi_{\overline{1}\overline{1}}(r^{-\frac{1}{2}}\theta^{\overline{1}})\circ(r^{-\frac{1}{2}}\theta^{\overline{1}})\right)\\\\
%&&+r^\frac{3}{2}\left(\Psi_{0\overline{1}}(r^{-1}\theta^0)\circ(r^{-\frac{1}{2}}\theta^{\overline{1}})+\Psi_{\overline{0}1}(r^{-1}\theta^{\overline{0}})\circ(r^{-\frac{1}{2}}\theta^1)\right)\\\\
%&&+O_e(r^2).
%\end{array}
%\end{displaymath} 
 
We have 
\begin{displaymath}
\begin{array}{rcl}
[e_0,e_{\overline{0}}]&=&\frac{1}{\sqrt{2}}\left(e_{\overline{0}}-e_0-r^{\frac{3}{2}}\Psi_{\overline{0}1}e_{\overline{1}}+r^\frac{3}{2}\Psi_{0\overline{1}}e_1\right)+O_e(r^2),\\

[e_0,e_1]&=&\frac{1}{\sqrt{2}}\left(\left(\frac{1}{2}-r\Phi_{1\overline{1}}\right)e_1-r\Phi_{11} e_{\overline{1}}-i\left(\tilde{\nabla}^\theta_{e_0}e_1-\tau(e_1)\right)\right)+O_e(r^2),\\

[e_0,e_{\overline{1}}]&=&\frac{1}{\sqrt{2}}\left(\left(\frac{1}{2}-r\Phi_{1\overline{1}}\right)e_{\overline{1}}-r\Phi_{\overline{1}\overline{1}} e_1-i\left(\tilde{\nabla}^\theta_{e_0}e_{\overline{1}}-\tau(e_{\overline{1}})\right)\right)+O_e(r^2),\\

[e_1, e_{\overline{1}}]&=&r^\frac{3}{2}\left(\Phi_{1\overline{1},\overline{1}}-\Phi_{\overline{1}\overline{1},1}\right)e_1-r^\frac{3}{2}\left(\Phi_{1\overline{1},1}-\Phi_{11,\overline{1}}\right)e_{\overline{1}}+O_e(r^2).
\end{array}
\end{displaymath}
Hence,
\begin{displaymath}
\begin{array}{rcl}
\nabla^g_{e_I}e_{\overline{I}}&=&2r\left(1-r\Phi_{1\overline{1}}\right)\partial_r\\
&&-r^2\left(\sqrt{2}\Psi_{0\overline{1}}+\Phi_{1\overline{1},\overline{1}}-\Phi_{\overline{1}\overline{1},1}\right)T_1\\
&&-r^2\left(\sqrt{2}\Psi_{\overline{0}1}+\Phi_{1\overline{1},1}-\Phi_{11,\overline{1}}\right)T_{\overline{1}}+O_T(r^\frac{5}{2}).
\end{array}
\end{displaymath}

We also have, for $\omega\in \Omega^1(X)\otimes\tilde{\varphi}^* TN$,
\begin{displaymath}
\begin{array}{rcl}
\nabla^{\tilde{\varphi}^*h}_{e_0}\omega(e_{\overline{0}})+\nabla^{\tilde{\varphi}^*h}_{e_{\overline{0}}}\omega(e_0)&=&r\omega(\partial_r)+r^2\nabla^{\tilde{\varphi}^*h}_{\partial_r}\omega(\partial_r)+r^2\nabla^{\tilde{\varphi}^*h}_R\omega(R)\\
&&-\sqrt{2}r^2\Psi_{0\overline{1}}\omega(T_1)-\sqrt{2}r^2\Psi_{\overline{0}1}\omega(T_{\overline{1}})+O_T(r^\frac{5}{2}),\\

\nabla^{\tilde{\varphi}^*h}_{e_1}\omega(e_{\overline{1}})&=&r\nabla^{\tilde{\varphi}^*h}_{T_1}\omega(T_{\overline{1}})-r^2\nabla^{\tilde{\varphi}^*h}_{T_1}\left(\Phi_{1\overline{1}}\omega(T_{\overline{1}})\right)-r^2\nabla^{\tilde{\varphi}^*h}_{T_1}\left(\Phi_{\overline{1}\overline{1}}\omega(T_1)\right)\\
&&-r^2\Phi_{1\overline{1}}\nabla^{\tilde{\varphi}^*h}_{T_1}\omega(T_{\overline{1}})-r^2\Phi_{11}\nabla^{\tilde{\varphi}^*h}_{T_{\overline{1}}}\omega(T_{\overline{1}})+O_T(r^3).

\end{array}
\end{displaymath}

\noindent The divergence is hence given by
\begin{displaymath}
\begin{array}{lll}
\delta^{g,h}\omega%&=&r(1-r\Phi_{1\overline{1}})\left(\omega(\partial_r)-\nabla^{\tilde{\varphi}^*h}_{T_1}\omega(T_{\overline{1}})-\nabla^{\tilde{\varphi}^*h}_{T_{\overline{1}}}\omega(T_1)\right)\\
%&&-r^2\nabla^{\tilde{\varphi}^*h}_{\partial_r}\omega(\partial_r)-r^2\nabla^{\tilde{\varphi}^*h}_R\omega(R)\\
%&&+2r^2\left(\nabla^{\tilde{\varphi}^* h}_{T_{\overline{1}}}\Phi_{11}\omega(T_{\overline{1}})+\nabla^{\tilde{\varphi}^* h}_{T_1}\Phi_{\overline{1}\overline{1}}\omega(T_1)\right)+O_T(r^\frac{5}{2})\\
&=&r(1-2r\Phi_{1\overline{1}})\left(\omega(\partial_r)-\nabla^{\tilde{\varphi}^*h}_{T_1}\omega(T_{\overline{1}})-\nabla^{\tilde{\varphi}^*h}_{T_{\overline{1}}}\omega(T_1)\right)-r^2\nabla^{\tilde{\varphi}^*h}_{\partial_r}\omega(\partial_r)-r^2\nabla^{\tilde{\varphi}^*h}_R\omega(R)\\
&&+4r^2\mathrm{Im}\left(\nabla^{\tilde{\varphi}^* h}_{T_{\overline{1}}}\left(\tau_1^{\overline{1}} \omega(T_{\overline{1}})\right)\right)+O_T(r^\frac{5}{2})\\
&=&{\delta^{g_0,h}\omega}-2r^2\Phi_{1\overline{1}}\nabla^{\tilde{\varphi}^*h}_{\partial_r}\delta^{g,h}\omega+4r^2\mathrm{Im}\left(\nabla^{\tilde{\varphi}^* h}_{T_{\overline{1}}}\left(\tau_1^{\overline{1}} \omega(T_{\overline{1}})\right)\right)+O_T(r^\frac{5}{2}).
\end{array}
\end{displaymath}

Then, by Theorem \ref{obs}, we have
\begin{displaymath}
\begin{array}{rcl}
P_1(\varphi)&=&\nabla^{\tilde{\varphi}^*h}_{\partial_r}\tilde{\delta}^{\theta,h}_bT\tilde{\varphi}|_{r=0}-\nabla^{\varphi^*h}_R R\varphi+4\mathrm{Im}\left(\nabla^{\varphi^* h}_{T_{\overline{1}}}\left(\tau_1^{\overline{1}} T_{\overline{1}}\right)\right)\varphi\\
&=&\delta^{\theta,h}_b\nabla^{\varphi^* h}\varphi_1+\mathcal{R}^h_{\varphi_1,T_1\varphi}T_{\overline{1}}\varphi+\mathcal{R}^h_{\varphi_1,T_{\overline{1}}\varphi}T_1\varphi-\nabla^{\varphi^*h}_R R\varphi+4\mathrm{Im}\left(\nabla^{\varphi^* h}_{T_{\overline{1}}}\left(\tau_1^{\overline{1}} T_{\overline{1}}\right)\right)\varphi,
\end{array}
\end{displaymath}
hence, since $\varphi_1=-\delta^{\theta,h}_b T\varphi$, the announced obstruction.
\end{proof}

Note that on functions, meaning that $N=\R$, $P_1$ reduces to a multiple of the CR Paneitz operator. Since the construction follows the ideas of Graham \emph{et al.}, this was expected. A similar phenomenon appears in the real case \cite{Ber13}.% Consequently, we have

%\begin{prop}
%If $n=1$, CR-pluriharmonic functions are CR-harmonic.
%\end{prop}
%\begin{proof}
%Indeed, CR-pluriharmonic functions are contained in the kernel of the CR Paneitz operator, see \cite{CCY15}.
%\end{proof}

\begin{ex}%\hfill
%\begin{itemize}
%\item[$\bullet$]
Let us consider $\mathrm{id}:(M,H,J,\theta)\rightarrow (M,g:=g_{J,\theta})$.\\ Since $\nabla^g_R R=0$ by Lemma 1.3. in \cite{DT06}, we have, using the Koszul formula:
\begin{displaymath}
\begin{array}{rcl}
\delta^{\theta,g}_bT\mathrm{id}&=&-\nabla^{g}_{T_1}T_{\overline{1}}-\nabla^{g}_{T_{\overline{1}}}T_1\\
&=&-g\left([T_1,T_{\overline{1}}],T_1\right)T_{\overline{1}}-g\left([T_{\overline{1}},T_1],T_{\overline{1}}\right)T_1\\
&&-g\left([R,T_{\overline{1}}],T_1\right)R-g\left([R,T_1],T_{\overline{1}}\right)R\\
&&-g\left([T_1,R],R\right)T_{\overline{1}}-g\left([T_{\overline{1}},R],R\right)T_1\\
&=&0,
\end{array}
\end{displaymath}
hence
 $$P_1(\mathrm{id})=4\mathrm{Im}\nabla^{g}_{T_{\overline{1}}}\left(\tau_{1}^{\overline{1}} T_{\overline{1}}\right).$$ Consequently, the identity is CR-harmonic if and only if $\mathrm{Im}\nabla^{g}_{T_{\overline{1}}}\left(\tau_{1}^{\overline{1}} T_{\overline{1}}\right)=0$. This is in particular verified when $\theta$ is normal, \emph{i.e.} when $\tau=0$.
%\item[$\bullet$] If $T_{\overline{1}}\varphi=0$, then, by Remark \ref{brackets}, $$\delta^{\theta,h}_b\nabla^{\varphi^*h}\delta^{\theta,h}_bT\varphi=-\nabla^{\varphi^*h}_R R\varphi+2i\nabla^{\varphi^*h}_{T_1} \left(\tau_{\overline{1}}^1 T_1\varphi\right)=-\nabla^{\varphi^*h}_R R\varphi+4\mathrm{Im}\left(\nabla^{\varphi^* h}_{T_{\overline{1}}}\left(\tau_1^{\overline{1}} T_{\overline{1}}\right)\right)\varphi,$$ hence $P_1(\varphi)=0$.
%\end{itemize}
\end{ex}

%As an example, let us consider $\mathrm{id}:M\rightarrow M$. We have $$P_1(\mathrm{id})=4\mathrm{Im}\left(\tau_{1,\overline{1}}^{\overline{1}} T_{\overline{1}}\right),$$ hence the identity is CR-harmonic if and only if $\tau_{1,\overline{1}}^{\overline{1}}=0$.

\section{Renormalized energy}\label{secren}

\subsection{Definition}

Let $\varphi\in C^\infty(M,N)$ and $\tilde{\varphi}$ be the extension of $\varphi$ constructed in Theorem \ref{obs}. For $\rho$ in $(0,\varepsilon)$, let $$E(\tilde{\varphi},\rho)=\frac{1}{2}\int_{(\rho,\varepsilon)\times M}\|T\tilde{\varphi}\|^2_{g,h} d\mathrm{vol}_g$$ be the energy of $\tilde{\varphi}$ in $(\rho,\varepsilon)\times M$. We have $$\|T\tilde{\varphi}\|_{g,h}^2=f_0 r + f_1 r^2+\ldots + f_nr^{n+1}+O(r^{n+2}\log r),$$ where $\forall k\leq n,$ $f_k$ depends only on $U_j$ for $j\leq k$ and on $g$ at order $k$ in $e^{(0)}$, and $$d\mathrm{vol}_g=r^{-n-2}\sqrt{\det g}dr\wedge\theta\wedge d\theta^n.$$
%, since $\tilde{\varphi}$ and $g$ have a Taylor development in powers of $r$ formally determined up to $O_T(r^n)$ and $O_e\left(r^{n+\frac{1}{2}}\right)$ respectively,
Consequently, $$\|T\tilde{\varphi}\|_{g,h}^2 d\mathrm{vol}_{g}=\left(a_0 r^{-n-1}+a_1 r^{-n} +\ldots + a_n r^{-1}+O(\log r)\right) dr\wedge\theta\wedge d\theta^n,$$ where $\forall k\leq n,$ $a_k$ depends only on $U_j$ for $j\leq k$ and on $g$ at order $k$. %We can write $$E(\tilde{\varphi},\rho)=\frac{1}{2}\int_{(\rho,\varepsilon)\times M}\left(e_0 r^{-n}+e_1 r^{1-n} +\ldots + e_{n}\right)dr\ \frac{d\mathrm{vol}_{g}}{r^n}+O(1),$$
Hence $E$ admits the development, when $\rho\rightarrow 0$, $$E(\tilde{\varphi},\rho)=E_0(\varphi) \rho^{-n}+E_1(\varphi) \rho^{1-n} +\ldots + E_{n-1}(\varphi)\rho^{-1} + F_n(\varphi)\log\rho + E_n(\varphi) + o(1),$$ where $\forall k\leq n-1,$ $E_k$ depends only on $U_j$ for $j\leq k$ and on $g$ at order $k$, and $F_n$ depends only on $U_j$ for $j\leq n$ and on $g$ at order $n$. The coefficient $F_n(\varphi)$ can be written as $$F_n(\varphi)=-\frac{1}{2}\int_M a_n \theta\wedge d\theta^n = -\frac{1}{2n!}\int_M \left.\partial_r^n\left(r^{n+1}\|T\tilde{\varphi}\|_{g,h}^2 d\mathrm{vol}_{g}\right)\right|_{r=0}.$$

By construction, $F_n$ is formally determined by $\varphi$, $g$ and $h$. Moreover, we have:

\begin{prop}\label{ene}
$F_n(\varphi)$ is a CR invariant: $$\hat{F}_n(\varphi)=F_n(\varphi).$$
\end{prop}

\begin{proof}
The proof is similar to the proof of Proposition \ref{crinv}. Indeed, if $\hat{r}=e^f r$, then
\begin{displaymath}
\begin{array}{rcl}
\|T\tilde{\varphi}\|_{g,h}^2 d\mathrm{vol}_{g}&=&\left(a_0 r^{-n-1}+a_1 r^{-n} +\ldots + a_n r^{-1}+a_{n+1}+O(r)\right) dr\wedge\theta\wedge d\theta^n\\
&=&\left(\hat{a}_0 \hat{r}^{-n-1}+\hat{a}_1 \hat{r}^{-n} +\ldots + \hat{a}_n\hat{r}^{-1}+\hat{a}_{n+1}+O(\hat{r})\right) d\hat{r}\wedge\theta\wedge d\theta^n,
\end{array}
\end{displaymath}
hence, when integrating over $(r=\rho,r=\varepsilon)\times M$,
\begin{displaymath}
\begin{array}{rcl}
E(\tilde{\varphi},\rho)&=&E_0(\varphi) \rho^{-n}+E_1(\varphi) \rho^{1-n} +\ldots + E_{n-1}(\varphi)\rho^{-1} + F_n(\varphi)\log\rho + E_n(\varphi) + o(1)\\
&=&\hat{E}_0(\varphi) \rho^{-n}+\hat{E}_1(\varphi) \rho^{1-n} +\ldots + \hat{E}_{n-1}(\varphi)\rho^{-1} + \hat{F}_n(\varphi)\log\rho + \hat{E}_n(\varphi) + o(1).
\end{array}
\end{displaymath}
Again, since the function $f$ itself has a Taylor expansion in $r$, all polynomial terms are mixed. However, the only $\log \rho$ term which appears when integrating with respect to $\hat{r}$ comes from the $\hat{r}^{-1}$ term. Hence the result.
\end{proof}

The principal term of $F_n(\varphi)$ is the following: since $$r^{n+1}\|T\tilde{\varphi}\|_{g,h}^2 d\mathrm{vol}_{g}=\left(\left<T_A\tilde{\varphi},T_{\overline{A}}\tilde{\varphi}\right>_h+r\|\partial_r\tilde{\varphi}\|_h^2\right) dr\wedge\theta\wedge d\theta^n+\text{lower order (in derivations of $\varphi$) terms},$$ we have
\begin{displaymath}
\begin{array}{rcl}
F_n(\varphi) &=& -\frac{1}{2n!}\int_M \left(\sum_{k=0}^n\binom{n}{k}\left<\delta^{\theta,h}_b\nabla^{\varphi^* h}\varphi_k,\varphi_{n-k}\right>_h+ n\sum_{k=0}^{n-1}\binom{n-1}{k}\left<\varphi_{k+1},\varphi_{n-k}\right>_h\right) \theta\wedge d\theta^n+\text{l.o.t.}\\
&=& \frac{(-1)^{n+1}}{2n!^2}\int_M \left<(\delta^{\theta,h}_b\nabla^{\varphi^*h})^{n-1}\delta^{\theta,h}_bT\varphi,\delta^{\theta,h}_bT\varphi\right>_h \theta\wedge d\theta^n+\text{lower order terms}.
\end{array}
\end{displaymath}

\begin{defi}
$F_n(\varphi)$ is called the \emph{renormalized energy} of $\varphi$.
\end{defi}

\begin{prop}
The gradient of $F_n(\varphi)$ is $\frac{1}{2n!}P_n(\varphi)$, that is to say, for all $\dot{\varphi}\in \Gamma(\varphi^*TN)$, $$d_\varphi F_n(\dot{\varphi})=\frac{1}{2n!}\int_M \left<\dot{\varphi},P_n(\varphi)\right>_h \theta\wedge d\theta^n.$$
\end{prop}

\begin{proof}
Let $\dot{\varphi}\in \Gamma(\varphi^*TN)$. Let $(\varphi_t)_{t\in[-1,1]}$ be a one-parameter family in $C^\infty(M,N)$ such that
\begin{displaymath}
\left\lbrace
\begin{array}{rcl}
\varphi_0=\varphi,\\
\partial_t\varphi_t|_{t=0}=\dot{\varphi}.
\end{array}
\right.
\end{displaymath}
Let us equip $X\times[-1,1]$ with the metric $\overline{g}=g+dt^2$ and let $\xi\in C^\infty(X\times[-1,1],N)$ be the map $$\forall p\in X,\ \forall t\in[-1,1],\quad \xi(p,t)=\tilde{\varphi}_t(p).$$

We then have 
\begin{displaymath}
\begin{array}{rcl}
\partial_t\|T\tilde{\varphi}_t\|^2_{g,h}&=&\partial_t\left(\|T\xi\|^2_{\overline{g},h}-\|\partial_t\xi\|^2_h\right)\\
&=&\left<\nabla^{\overline{g},h}_{\partial_t}T\xi,T\xi\right>_{\overline{g},h}-\left<\nabla^{\xi^*h}_{\partial_t}\partial_t\xi,\partial_t\xi\right>_h\\
&=&\left<\nabla^{\xi^*h}_{\partial_t}e_I\xi,e_{\overline{I}}\xi\right>_h\\
&=&\left<\nabla^{\xi^*h}_{e_I}\partial_t\xi,e_{\overline{I}}\xi\right>_h\\
%&=& \left<\nabla^{\tilde{\varphi}_t^*h}_{\partial_t}e_I\tilde{\varphi}_t,e_{\overline{I}}\tilde{\varphi}_t\right>_h\\
&=& \left<\nabla^{\tilde{\varphi}_t^*h}_{e_I}\partial_t\tilde{\varphi}_t,e_{\overline{I}}\tilde{\varphi}_t\right>_h\\
&=&e_I\left<\partial_t\tilde{\varphi}_t,e_{\overline{I}}\tilde{\varphi}_t\right>_h-\left<\partial_t\tilde{\varphi}_t,\nabla^{\tilde{\varphi}_t^*h}_{e_I}e_{\overline{I}}\tilde{\varphi}_t\right>_h,
\end{array}
\end{displaymath}
hence
\begin{displaymath}
\begin{array}{rcl}
\partial_t E(\tilde{\varphi}_t,\rho)|_{t=0}&=&\frac{1}{2}\int_{(\rho,\varepsilon)\times M}\left(e_I\left<\partial_t\tilde{\varphi}_t|_{t=0},e_{\overline{I}}\tilde{\varphi}\right>_h-\left<\partial_t\tilde{\varphi}_t|_{t=0},\nabla^{\tilde{\varphi}^*h}_{e_I}e_{\overline{I}}\tilde{\varphi}\right>_h\right) d\mathrm{vol}_g.
\end{array}
\end{displaymath}

There is no $\log\rho$ term in the second part, and
\begin{displaymath}
\begin{array}{rcl}
\frac{1}{2}\int_{(\rho,\varepsilon)\times M}e_I\left<\partial_t\tilde{\varphi}_t|_{t=0},e_{\overline{I}}\tilde{\varphi}\right>_h d\mathrm{vol}_g&=&\frac{1}{2}\int_M \rho^{-n}\left<\partial_t\tilde{\varphi}_t|_{t=0},\partial_\rho\tilde{\varphi}\right>_h \theta\wedge d\theta^n + \text{lower order terms},
\end{array}
\end{displaymath}
whose $\log\rho$ term is $$\frac{1}{2n!}\int_M \left<\dot{\varphi},P_n(\varphi)\right>_h \theta\wedge d\theta^n,$$
hence the result.
\end{proof}

\subsection{Explicit energy in dimension $3$}
\label{exen}

Here again, when $n=1$, \emph{i.e.} $\dim(M)=3$, knowing the asymptotic development of $g$ at order $\frac{3}{2}$ in $e^{(0)}$ allows for an explicit computation of the renormalized energy.

\begin{theo}\label{corpsB2}
We have $$F_1(\varphi)=-\frac{1}{2}\int_M \left(\|\delta^{\theta,h}_bT\varphi\|^2_h+\|R\varphi\|^2_h-4\mathrm{Im}\left(\tau_1^{\overline{1}} \|T_{\overline{1}}\varphi\|^2_h\right)\right)\theta\wedge d\theta.$$
\end{theo}

\begin{proof}

We have
\begin{displaymath}
\begin{array}{rcl}
\|T\tilde{\varphi}\|^2_{g,h} &=& 2\left<e_0\tilde{\varphi},e_{\overline{0}}\tilde{\varphi}\right>_h+2\left<e_1\tilde{\varphi},e_{\overline{1}}\tilde{\varphi}\right>_h\\
&=&2r\left<T_1\varphi,T_{\overline{1}}\varphi\right>_h\\
&&+r^2\left(\|\varphi_1\|^2_h+\|R\varphi\|^2_h-4\Phi_{1\overline{1}}\left<T_1\varphi,T_{\overline{1}}\varphi\right>_h-2\Phi_{11}\|T_{\overline{1}}\varphi\|^2_h-2\Phi_{\overline{1}\overline{1}}\|T_1\varphi\|^2_h\right)\\
&&+O(r^\frac{5}{2}),
\end{array}
\end{displaymath}
and $$d\mathrm{vol}_g=\left(1+2r\Phi_{1\overline{1}}+O(r^{2})\right)r^{-3}dr\wedge\theta\wedge d\theta.$$
Consequently,
\begin{displaymath}
\begin{array}{lll}
r^2\|T\tilde{\varphi}\|^2_g d\mathrm{vol}_g &=&\left(2\left<T_1\varphi,T_{\overline{1}}\varphi\right>_h\right.\\
&&+r\left(\|\varphi_1\|^2_h+\|R\varphi\|^2_h-2\Phi_{11}\|T_{\overline{1}}\varphi\|^2_h-2\Phi_{\overline{1}\overline{1}}\|T_1\varphi\|^2_h\right)\\
&&\left.+O(r^{\frac{3}{2}})\right)dr\wedge\theta\wedge d\theta,
\end{array}
\end{displaymath}
and finally
\begin{displaymath}
\begin{array}{rcl}F_1(\varphi)&=&-\frac{1}{2}\int_M \left(\|\varphi_1\|^2_h+\|R\varphi\|^2_h-4\mathrm{Im}\left(\tau_1^{\overline{1}} \|T_{\overline{1}}\varphi\|^2_h\right)\right)\theta\wedge d\theta.
%&=\frac{1}{2}\int_M &\left(-\|\delta^{W_0,h}T\varphi\|^2_h+\|R\varphi\|^2_h+8iA_{11}\left(\|T_1\varphi\|^2_h+\|T_{\overline{1}}\varphi\|^2_h\right)\right)\theta\wedge d\theta.
\end{array}
\end{displaymath}

\end{proof}

As an example, for $\mathrm{id}:(M,H,J,\theta)\rightarrow (M,g_{J,\theta})$, we have $$F_1(\mathrm{id})=-\frac{1}{2}\mathrm{Vol}(M,\theta).$$

\section{Further computations in the general case}\label{secfur}

We give here a more precise computation for $\delta^{g,h}\omega$ and $r^{n+1}\|T\tilde{\varphi}\|^2_g d\mathrm{vol}_g$ in the general case, using Theorem \ref{poinc}. We show that this computation does not allow for an explicit expression of the obstruction and of the renormalized energy respectively.

\subsection{Computation of the divergence}

%Let us consider the development
%
%$$g=g_0+\Phi_{AB}\theta^A\circ\theta^B+O_e(r^\frac{3}{2}).$$

By Theorem \ref{poinc}, we have $$g=g_0+\Phi_{AB}\theta^A\circ\theta^B+O_e(r^\frac{3}{2}),$$ where, denoting by $R_{\alpha\overline{\beta}}$ the components of the Webster Ricci tensor, $$\Phi_{\alpha\overline{\beta}}=-\frac{1}{n+2}\left(R_{\alpha\overline{\beta}}-\frac{\mathrm{Scal}_W(J,\theta)}{2(n+1)}\delta_{\alpha\overline{\beta}}\right),\quad\text{and}\quad\Phi_{\alpha\beta}=-i \tau_\alpha^{\overline{\beta}}.$$% where 

\noindent By Proposition \ref{comp2}, we can equip $\lbrace r\rbrace\times H$ with a complex structure $J_r=J_0+rJ_1+O_T(r^2)$, with $$J_1 T_\alpha=-2\Phi_{\alpha\beta}T_{\overline{\beta}}.$$ An orthonormal basis of $T^{1,0}X$ with respect to $g$ %, provided $Jr\partial_r:=rR$,
is given by
\begin{displaymath}
\begin{array}{lll}
(e_0,e_\alpha)&:=&\left(r\partial_r-i rR, \left(\delta_{\alpha\overline{\beta}}-r\Phi_{\alpha\overline{\beta}}\right) r^\frac{1}{2}T_\beta-r\Phi_{\alpha\beta}r^\frac{1}{2} T_{\overline{\beta}}\right)+O_e(r^\frac{3}{2}).
\end{array}
\end{displaymath}
Now, $g$ can be rewritten as $$g=(r^{-1}\theta^0)\circ(r^{-1}\theta^{\overline{0}})+(r^{-\frac{1}{2}}\theta^\alpha)\circ(r^{-\frac{1}{2}}\theta^{\overline{\alpha}})+r\Phi_{AB}(r^{-\frac{1}{2}}\theta^A)\circ(r^{-\frac{1}{2}}\theta^B)+O_e(r^\frac{3}{2}).$$

%Let us denote by $\rho$ the matrix of $T^W(R,\cdot)$ in the basis $h$: $$T^W(R,h_\alpha)=\rho_{\alpha\beta}h_\beta+\rho_{\alpha\tilde{\beta}}Jh_\beta,$$ $$T^W(R,Jh_\alpha)=\rho_{\tilde{\alpha}\beta}h_\beta+\rho_{\tilde{\alpha}\tilde{\beta}}Jh_\beta.$$ Since $T^W(R,\cdot)$ and $J$ anticommute, $\rho$ is block-symmetric and -trace-free, \emph{i.e.} $$\rho_{\tilde{\alpha}\beta}=\rho_{\alpha\tilde{\beta}}\quad\mathrm{and}\quad\rho_{\tilde{\alpha}\tilde{\beta}}=-\rho_{\alpha\beta}.$$
We have, modulo $O_e(r^\frac{3}{2})$,
\begin{displaymath}
\begin{array}{rcl}
[e_0,e_{\overline{0}}]&=&\frac{1}{\sqrt{2}}\left(e_{\overline{0}}-e_0\right),\\

[e_0,e_\alpha]&=&\frac{1}{\sqrt{2}}\left(\frac{1}{2}e_\alpha-r\Phi_{\alpha\overline{\beta}}e_\beta-r\Phi_{\alpha\beta} e_{\overline{\beta}}-i\left(\tilde{\nabla}^\theta_{e_0}e_\alpha-\tau(e_\alpha)\right)\right),\\

[e_0,e_{\overline{\alpha}}]&=&\frac{1}{\sqrt{2}}\left(\frac{1}{2}e_{\overline{\alpha}}-r\Phi_{\overline{\alpha}\beta}e_{\overline{\beta}}-r\Phi_{\overline{\alpha}\overline{\beta}} e_\beta-i\left(\tilde{\nabla}^\theta_{e_0}e_{\overline{\alpha}}-\tau(e_{\overline{\alpha}})\right)\right),\\

[e_\alpha,e_\beta]&=&r^\frac{3}{2}\left(\Phi_{\alpha\overline{\delta},\beta}-\Phi_{\beta\overline{\delta},\alpha}\right)e_\delta+r^\frac{3}{2}\left(\Phi_{\alpha\delta,\beta}-\Phi_{\beta\delta,\alpha}\right)e_{\overline{\delta}},\\

[e_\alpha,e_{\overline{\beta}}]&=&r^\frac{3}{2}\left(\Phi_{\alpha\overline{\delta},\overline{\beta}}-\Phi_{\overline{\beta}\overline{\delta},\alpha}\right)e_\delta-r^\frac{3}{2}\left(\Phi_{\delta\overline{\beta},\alpha}-\Phi_{\alpha\delta,\overline{\beta}}\right)e_{\overline{\delta}}.

\end{array}
\end{displaymath}
Hence,
\begin{displaymath}
\begin{array}{ll}

\nabla^g_{e_i}e_{\overline{i}}+\nabla^g_{e_{\overline{i}}}e_i=&r\left(n+1-2r\Phi_{\alpha\overline{\alpha}}\right)\partial_r\\
&-r^2\left(2\Phi_{\beta\overline{\beta},\overline{\alpha}}-\Phi_{\overline{\alpha}\overline{\beta},\beta}-\Phi_{\overline{\alpha}\beta,\overline{\beta}}\right)T_\alpha\\
&-r^2\left(2\Phi_{\beta\overline{\beta},\alpha}-\Phi_{\alpha\beta,\overline{\beta}}-\Phi_{\alpha\overline{\beta},\beta}\right)T_{\overline{\alpha}}+O_T(r^\frac{5}{2}).
\end{array}
\end{displaymath}
with $\Phi_{\alpha\overline{\alpha}}=- \mathrm{tr}(\mathrm{S_\theta})=-\frac{\mathrm{Scal}_W}{2(n+1)}$. %and $\rho_{\alpha\alpha}+\rho_{\tilde{\alpha}\tilde{\alpha}}=0$ by the consideration above.

Also,
\begin{displaymath}
\begin{array}{rcl}
\nabla^{\tilde{\varphi}^*h}_{e_0}\omega(e_{\overline{0}})+\nabla^{\tilde{\varphi}^*h}_{e_{\overline{0}}}\omega(e_0)&=&r\omega(\partial_r)+r^2\nabla^{\tilde{\varphi}^*h}_{\partial_r}\omega(\partial_r)+r^2\nabla^{\tilde{\varphi}^*h}_R\omega(R)+O_T(r^2),\\

\nabla^{\tilde{\varphi}^*h}_{e_\alpha}\omega(e_{\overline{\alpha}})&=&r\nabla^{\tilde{\varphi}^*h}_{T_\alpha}\omega(T_{\overline{\alpha}})-r^2\nabla^{\tilde{\varphi}^*h}_{T_\alpha}\left(\Phi_{\overline{\alpha}\beta}\omega(T_{\overline{\beta}})\right)-r^2\nabla^{\tilde{\varphi}^*h}_{T_\alpha}\left(\Phi_{\overline{\alpha}\overline{\beta}}\omega(T_\beta)\right)\\
&&-r^2\Phi_{\alpha\overline{\beta}}\nabla^{\tilde{\varphi}^*h}_{T_\beta}\omega(T_{\overline{\alpha}})-r^2\Phi_{\alpha\beta}\nabla^{\tilde{\varphi}^*h}_{T_{\overline{\beta}}}\omega(T_{\overline{\alpha}})+O_T(r^{\frac{5}{2}}).
\end{array}
\end{displaymath}

Coming back to the divergence, we have
\begin{displaymath}
\begin{array}{lll}
\delta^{g,h}\omega&=&r(n-2r\Phi_{\alpha\overline{\alpha}})\omega(\partial_r)-r^2\nabla^{\tilde{\varphi}^*h}_{\partial_r}\omega(\partial_r)-r^2\nabla^{\tilde{\varphi}^*h}_R(\omega(R))\\
&&-r(1-2r\Phi_{\alpha\overline{\alpha}})\left(\nabla^{\tilde{\varphi}^*h}_{T_\alpha}\omega(T_{\overline{\alpha}})+\nabla^{\tilde{\varphi}^*h}_{T_{\overline{\alpha}}}\omega(T_\alpha)\right)\\
&&+2r^2\left(\nabla^{\tilde{\varphi}^*h}_{T_{\overline{\beta}}}\left(\Phi_{\alpha\beta}\omega(T_{\overline{\alpha}})\right)+\nabla^{\tilde{\varphi}^*h}_{T_\beta}\left(\Phi_{\overline{\alpha}\overline{\beta}}\omega(T_\alpha)\right)\right)\\
&&+2r^2\left(\nabla^{\tilde{\varphi}^*h}_{T_\alpha}\Phi_{\overline{\alpha}\beta}\omega(T_{\overline{\beta}})+\nabla^{\tilde{\varphi}^*h}_{T_{\overline{\alpha}}}\Phi_{\alpha\overline{\beta}}\omega(T_\beta)-\nabla^{\tilde{\varphi}^*h}_{T_\alpha}\Phi_{\beta\overline{\beta}}\omega(T_{\overline{\alpha}})-\nabla^{\tilde{\varphi}^*h}_{T_{\overline{\alpha}}}\Phi_{\beta\overline{\beta}}\omega(T_\alpha)
\right)+O_T(r^2).
\end{array}
\end{displaymath}

The term of order $2$ is consequently not known, which does not allow for an explicit computation of $P_n$. Note that $$\nabla^{\tilde{\varphi}^*h}_{T_{\overline{\beta}}}\left(\Phi_{\alpha\beta}\omega(T_{\overline{\alpha}})\right)+\nabla^{\tilde{\varphi}^*h}_{T_\beta}\left(\Phi_{\overline{\alpha}\overline{\beta}}\omega(T_\alpha)\right)=2\mathrm{Im}\left(\nabla^{\tilde{\varphi}^*h}_{T_{\overline{\beta}}}\left(\tau_\alpha^{\overline{\beta}}\omega(T_{\overline{\alpha}})\right)\right),$$
and that the potentially hidden $r^2$ terms are necessarily of the form $C^\alpha r^2\omega(T_\alpha)+D^{\overline{\alpha}}r^2\omega(T_{\overline{\alpha}})$.

\subsection{Computation of the integrand of the energy}

We have
\begin{displaymath}
\begin{array}{rcl}
\|T\tilde{\varphi}\|^2_{g,h} &=& 2\left<e_0\tilde{\varphi},e_{\overline{0}}\tilde{\varphi}\right>_h+2\left<e_\alpha\tilde{\varphi},e_{\overline{\alpha}}\tilde{\varphi}\right>_h\\
&=&2r\left<T_\alpha\varphi,T_{\overline{\alpha}}\varphi\right>_h\\
&&+r^2\left(\|\varphi_1\|^2_h+\|R\varphi\|^2_h-2\Phi_{\alpha\beta}\left<T_{\overline{\alpha}}\varphi,T_{\overline{\beta}}\varphi\right>_h-2\Phi_{\overline{\alpha}\overline{\beta}}\left<T_\alpha\varphi,T_\beta\varphi\right>_h\right.\\
&&\left.-2\Phi_{\alpha\overline{\beta}}\left<T_{\overline{\alpha}}\varphi,T_\beta\varphi\right>_h-2\Phi_{\overline{\alpha}\beta}\left<T_\alpha\varphi,T_{\overline{\beta}}\varphi\right>_h\right)\\
&&+O(r^2),
\end{array}
\end{displaymath}
and $$d\mathrm{vol}_g=\left(1+2r\Phi_{\alpha\overline{\alpha}}+O(r^\frac{3}{2})\right)r^{-n-2}dr\wedge\theta\wedge d\theta^n.$$

Consequently,

\begin{displaymath}
\begin{array}{rcl}
r^{n+1}\|T\tilde{\varphi}\|^2_g d\mathrm{vol}_g &=&\left(2\left<T_\alpha\varphi,T_{\overline{\alpha}}\varphi\right>_h\right.\\
&&+r\left(\|\varphi_1\|^2_h+\|R\varphi\|^2_h-2\Phi_{\alpha\beta}\left<T_{\overline{\alpha}}\varphi,T_{\overline{\beta}}\varphi\right>_h-2\Phi_{\overline{\alpha}\overline{\beta}}\left<T_\alpha\varphi,T_\beta\varphi\right>_h\right.\\
&&\left.-2\Phi_{\alpha\overline{\beta}}\left<T_{\overline{\alpha}}\varphi,T_\beta\varphi\right>_h-2\Phi_{\overline{\alpha}\beta}\left<T_\alpha\varphi,T_{\overline{\beta}}\varphi\right>_h+4\Phi_{\alpha\overline{\alpha}}\left<T_\alpha\varphi,T_{\overline{\alpha}}\varphi\right>_h\right)\\
&&\left.+O(r)\right)dr\wedge\theta\wedge d\theta^n.
\end{array}
\end{displaymath}

The term of order $1$ is consequently not known, which does not allow for an explicit computation of $F_n$.

\section{Relation with the Fefferman bundle in dimension $3$}\label{secfef}

We describe here the correspondance between the obstruction to CR-harmonicity on a given CR $3$-manifold and the obstruction to conformal-harmonicity on its Fefferman bundle. It generalizes the Appendix B. of \cite{CY13}.

%We extend the results of Bérard to the Lorentzian setting . 
Let $(M,H,J)$ be a compact strictly pseudoconvex CR $3$-manifold and let $(N,h)$ be a Riemannian manifold. Let $(F,g_F)$ be the Fefferman bundle of $(M,H,J)$. For a detailed construction of the Fefferman bundle, see \cite{Far86,Lee86,Her09}. Let $\pi:F\rightarrow M$ be the natural bundle projection. Let $\theta$ be a positive contact form on $M$ and let $\varpi$ be the $S^1$-invariant connection $1$-form induced by the Weyl structure attached to $\theta$ on $F$. The \emph{Fefferman metric} attached to $\theta$ on $F$ is the Lorentzian metric $$g_F=i\varpi\circ\pi^*\theta+\frac{1}{2}\pi^*\gamma.$$% where $\lambda\circ\mu:=\lambda\otimes\mu+\mu\otimes\lambda$.

By analogy with the Riemannian case \cite{Ber13}, given $\varphi\in C^\infty(F,N)$, the obstruction to the existence of a smooth harmonic extension of $\varphi$ on the interior of $(F,g_F)$ is given by $$P_F(\varphi)=-\frac{1}{16}\left(\delta^{g_F,h}\nabla^{\varphi^*h}\delta^{g_F,h}T\varphi-\delta^{g_F,h}\left(2\mathrm{Ric}_{g_F}-\frac{2}{3}\mathrm{Scal}_{g_F}\right)T\varphi+S(\delta^{g_F,h}T\varphi)\right),$$
where $\mathrm{Ric}_{g_F}$ is understood as an endomorphism of $TF$, and $\mathrm{Ric}_{g_F}T\varphi:=T\varphi\left(\mathrm{Ric}_{g_F}(\cdot)\right)$, and $$S(X):=\sum_{i=1}^4 \mathcal{R}^h_{X,T\varphi(e_i)}T\varphi(e_i).$$

\begin{prop}
For all $\varphi\in C^\infty(M^3,N)$,
$$\pi_*\left(\delta^{g_F,h}\nabla^{\varphi^*h}\delta^{g_F,h}T(\pi^*\varphi)\right)=4\delta^{\theta,h}_b \nabla^{\varphi^*h}\delta^{\theta,h}_b T\varphi,$$
$$\pi_*\left(\delta^{g_F,h}\left(2\mathrm{Ric}_{g_F}-\frac{2}{3}\mathrm{Scal}_{g_F}\right)T(\pi^*\varphi)\right)=-4\nabla^{\varphi^* h}_R R\varphi+16\mathrm{Im}\left(\nabla^{\varphi^* h}_{T_{\overline{1}}}\left(\tau_1^{\overline{1}} T_{\overline{1}}\right)\right)\varphi,$$
and for $X$ in $TN$,
$$\pi_*\left(S((\pi^*\varphi)^* X)\right)=4S_b(\varphi^* X).$$
\end{prop}

\begin{proof}
The first and third equalities are straightforward from the expression of $g_F$. The second equality comes from the fact that, see \cite{Lee86}, $$\mathrm{Sch}_{g_F}=%-\frac{\mathrm{Scal}_W}{2}i\varpi\circ\theta-2(\varpi^2+K\theta^2)+\left(\mathrm{Sch}_W+\frac{\mathrm{Scal}_W}{4}\gamma\right)-\gamma(J\tau \cdot,\cdot)
-\varpi^2-\bold{S}\theta^2+\frac{1}{2}\mathrm{Sch}_W-\frac{1}{2}\gamma(J\tau \cdot,\cdot)+\frac{1}{2}\bold{T}J\circ\theta,$$
where $$\bold{T}=\frac{1}{3}\left(\frac{1}{4}d_b\mathrm{Scal}_W+i\delta\tau\right) \quad \mathrm{and}\quad \bold{S}=\delta\bold{T}-|\mathrm{Sch}_W|^2+|\tau|^2.$$

Indeed, since $\mathrm{Scal}_{g_F}=3\mathrm{Scal}_W\ $  and  $\ \mathrm{Sch}_W=\frac{1}{4}\mathrm{Scal}_W\gamma$, we have then
\begin{displaymath}
\begin{array}{rcl}
2\mathrm{Ric}_{g_F}-\frac{2}{3}\mathrm{Scal}_{g_F}g_F &=& 4\mathrm{Sch}_{g_F}-\frac{1}{3}\mathrm{Scal}_{g_F} g_F\\
&=& -4\varpi^2-4\bold{S}\theta^2-2\gamma(J\tau \cdot,\cdot)+2\bold{T}J\circ\theta-\mathrm{Scal}_W i\varpi\circ\theta,
\end{array}
\end{displaymath}
which gives the second equality.
\end{proof}

From the latter comes directly the

\begin{theo}
For all $\varphi\in C^\infty(M^3,N)$,
$$\pi_*\left(P_F(\pi^*\varphi)\right)=\frac{1}{4} P_1(\varphi).$$
\end{theo}
\noindent In particular, a map $\varphi:M\rightarrow N$ is CR-harmonic if and only if $\pi^*\varphi$ is conformal-harmonic.

\bibliographystyle{alpha}

\bibliography{../Biblio}

\end{document}